\documentclass[12pt,reqno]{amsart}
\usepackage{amsmath}
\usepackage{amssymb}
\usepackage{verbatim}
\usepackage{mathrsfs}
\usepackage{graphicx} 
\usepackage{caption}
\usepackage{subcaption}
\usepackage{hyperref}
\usepackage{color}
\usepackage[noadjust]{cite}

\usepackage[top=0.5in,bottom=0.5in,left=0.7in,right=0.7in]{geometry}

\newtheorem{lemma}{Lemma}

\newtheorem{cor}[lemma]{Corollary}
\newtheorem{theorem}[lemma]{Theorem}

\newtheorem{exmp}{Example}

\newtheorem{definition}[lemma]{Definition}

\numberwithin{equation}{section}
\numberwithin{lemma}{section}

\newcommand{\C}{\mathbb{C}}    
\newcommand{\N}{\mathbb{N}}    
\newcommand{\PL}{\mathbb{P}}   
\newcommand{\R}{\mathbb{R}}    
\newcommand{\Z}{\mathbb{Z}}    

\newcommand{\wh}{\widehat}
\renewcommand{\le}{\leqslant}
\renewcommand{\ge}{\geqslant}
\newcommand{\bs}{\backslash}
\newcommand{\ol}{\overline}
\newcommand{\la}{\langle}
\newcommand{\ra}{\rangle}
\newcommand{\bo}{\mathcal{O}} 

\newcommand{\pp}{\mathsf{p}}

\newcommand{\DG}{{\mathsf{Diag}}}

\newcommand{\FF}{\mathsf{F}}

\newcommand{\er}{\eqref}

\newcommand{\mra}{{\mathring{a}}}
\newcommand{\mrb}{{\mathring{b}}}
\newcommand{\mrc}{{\mathring{c}}}
\newcommand{\mrd}{{\mathring{d}}}

\newcommand{\mrv}{{\mathring{v}}}

\newcommand{\mrta}{{\mathring{\tilde{a}}}}
\newcommand{\mrtb}{{\mathring{\tilde{b}}}}

\newcommand{\mrphi}{\mathring{\phi}}
\newcommand{\mrtvgu}{\mathring{\tilde{\vgu}}}
\newcommand{\bpphi}{\breve{\phi}}
\newcommand{\mrtphi}{\mathring{\tilde{\phi}}}
\newcommand{\mrvgu}{\mathring{\vgu}}

\newcommand{\mrgam}{\mathring{\gamma}}

\newcommand{\tpsi}{\tilde{\psi}}
\newcommand{\tphi}{\tilde{\phi}}

\newcommand{\tvgu}{\tilde{\vgu}}

\newcommand{\bpa}{\breve{a}}
\newcommand{\bpb}{\breve{b}}

\newcommand{\bpta}{\breve{\tilde{a}}}
\newcommand{\bptb}{\breve{\tilde{b}}}

\newcommand{\bpvgu}{\breve{\vgu}}
\newcommand{\bptvgu}{\breve{\tilde{\vgu}}}
\newcommand{\bptphi}{\breve{\tilde{\phi}}}

\newcommand{\hta}{\widehat{\tilde{a}}}
\newcommand{\htb}{\widehat{\tilde{b}}}

\newcommand{\httheta}{\widehat{\tilde{\theta}}}

\newcommand{\bp}{ \begin{proof} }
	\newcommand{\ep}{\hfill \end{proof} }
\newcommand{\be}{ \begin{equation} }
\newcommand{\ee}{ \end{equation} }
\newcommand{\tp}{\mathsf{T}}
\newcommand{\dm}{\mathsf{M}} 
\newcommand{\dn}{\mathsf{N}} 
\newcommand{\vgu}{\upsilon} 
\newcommand{\Vgu}{\Upsilon} 

\newcommand{\sr}{\operatorname{sr}}  
\newcommand{\vmo}{\operatorname{vm}}
\newcommand{\bvmo}{\operatorname{bvm}}

\newcommand{\td}{\pmb{\delta}}  

\newcommand{\sd}{\mathcal{S}}  
\newcommand{\tz}{\mathcal{T}}  
\newcommand{\cM}{\mathcal{M}}

\newcommand{\cN}{\mathcal{N}}
\newcommand{\bpo}{\operatorname{bo}} 

\newcommand{\om}[1]{\omega_{#1}}
\newcommand{\ga}[1]{\gamma_{#1}}
\newcommand{\ka}[1]{\mathring{\gamma}_{#1}}

\newcommand{\dN}{\mathbb{N}^d}

\newcommand{\dR}{\mathbb{R}^d}

\newcommand{\dZ}{\mathbb{Z}^d}

\newcommand{\dLp}[1]{L_{#1}(\mathbb{R}^d)}

\newcommand{\dLrs}[3]{(L_{#1}(\mathbb{R}^d))^{#2\times #3}}

\newcommand{\vertiii}[1]{{\left\vert\kern-0.25ex\left\vert\kern-0.25ex\left\vert #1 
		\right\vert\kern-0.25ex\right\vert\kern-0.25ex\right\vert}}

\newcommand{\dlp}[1]{l_{#1}(\mathbb{Z}^d)}

\newcommand{\dlrs}[3]{(l_{#1}(\mathbb{Z}^d))^{#2\times #3}}
\newcommand{\dsq}{l(\mathbb{Z}^d)} 

\begin{document}
\title[Balanced Dual Framelets]
{A structural characterization of Compactly Supported OEP-based balanced dual multiframelets}

\author{Ran Lu}
\address{College of Science, Hohai University, Nanjing, China 211100.
\quad  {\tt rlu3@hhu.edu.cn}}

\thanks{The research of the author was supported by the National Natural Science Foundation of China under grants 12201178 and 12271140.
}

\makeatletter \@addtoreset{equation}{section} \makeatother
\begin{abstract}

	Compared to scalar framelets, multiframelets have certain advantages, such as relatively smaller supports on generators, high vanishing moments, etc.  The balancing property of multiframelets is very desired, as it reflects how efficient vector-valued data can be processed under the corresponding discrete multiframelet transform. Most of the literature studying balanced multiframelets is from the point of view of the function setting, but very few approaches are from the aspect of multiframelet filter banks. In this paper, we study structural characterizations of balanced dual multiframelets from the point of view of the Oblique Extension Principle (OEP). The OEP naturally connects framelets with filter banks, which makes it a very handy tool for analyzing the properties of framelets. With the OEP, we shall characterize compactly supported balanced dual multiframemets through the concept of balanced moment correction filters, which is the key notion that will be introduced in our investigation. The results of this paper demonstrate what essential structures a balanced dual multiframelet has in the most general setting, and bring us a more complete picture to understand balanced multiframelets and their underlying discrete multiframelet transforms.

\end{abstract}


 \maketitle

\pagenumbering{arabic}


\section{Introduction}
\label{sec:intro}

Multiwavelets and multiframelets are of interest in applications such as image process and numerical algorithms. Compared to scalar wavelets and framelets, multiwavelets and multiframelets have relatively smaller support on the generators and high vanishing moments. For literature on multiwavelets and multiframelets, we refer the readers to \cite{chs02,cj00,cj03,dhrs03,han97,han03,han09,han10,hl19pp,hl20pp,hm03,lv98,mo06,sel00,sel01} and many references therein. Quite often, multiwavelets and multiframelets are derived from refinable vector functions through extension principles. The most popular choices are the unitary extension principle (UEP) (\cite{rs97}) and its generalization oblique extension principle (OEP) (\cite{dhrs03,hm03,hl20pp}). With extension principles, wavelets and framelets can be characterized under either the function setting or the discrete filter bank setting. However, when it comes to multiwavelets and multiframelets, there is some discrepancy in the approximation properties under the two aforementioned settings, and this issue is known as the balancing property in the literature. In this paper, we will perform a comprehensive investigation on characterizing the structures of balanced multiframelets in the most general setting.

There has been a growing interest in investigating framelets and multiframelets in multi-dimensions these years, and several related progress has been made. For example, quasi-tight framelets with high order vanishing moments or directionality derived from scalar refinable functions through the unitary extension principle have been studied in \cite{dh18pp}; methods of constructing symmetric quincunx tight-framelets have been given in \cite{hjsz18}; multivariate quasi-tight multiframelts with high balancing orders has been investigated in \cite{hl20pp}; multivariate symmetric interpolatory multiframletes have been studied in \cite{k22}. Multiramelets generalize multiwavelets by adding redundancy to the system, and this offers much more flexibility for their construction. To process vector-valued data efficiently by a discrete multiwavelets/multiframelet transform (DmWT/DmFT), it is desired that the multiwavelets/multiframelet is balanced. For comprehensive motivations and backgrounds on the balancing property of multiwavelets, see \cite{cj00,cj03,hr98,lv98, sel98,sel00,xhgs96} and many references therein for more details. Although the balancing property of multiwavelets has been extensively investigated in the literature, the corresponding theory for multiframelets is still far from being well developed. We are only aware of the following work on the balancing property of multiframelets: \cite{han09} studied univariate balanced dual multiframelets and provided a rough guideline for construction; \cite{han10} investigated structural properties of multivariate balanced dual multiframelets; \cite{hl19pp,hl20pp} studied the existence and properties of univariate and multivariate balanced quasi-tight multiframelets respectively. Moreover, most of the approaches studying the balancing property of multiwavelets are from the point of view of functions, but it seems more natural to study the balancing property via the aspect of discrete filter banks (see e.g. \cite{han09,han10,hl19pp,hl20pp}), as the balancing property directly reflects the efficiency and the approximation property of a discrete multiframelet transform. In this paper, we shall follow the filter bank approach and provide a simple structural characterization of the balancing property of multiframelets in arbitrary dimensions. Our study is based on the popular oblique extension principle (OEP), which demonstrates the intrinsic connections of multiframelets systems and their corresponding multiframelet filter banks, and provides an elegant guideline to construct multiframelets with high vanishing moments and high balancing orders.

The structure of this paper is organized as the following: In Section~\ref{sec:pre}, we provide some necessary preliminaries for our later investigation, including: a brief review on dual framelets; the balancing property of a discrete multiframelet transforms; some basic theory on a recetly developed normal form of a matrix-valued filter which greatly facilitates our study of multiframelets. In Section~\ref{sec:cha}, we present our first main result Theorem~\ref{thm:dfrt}, which states that a balanced dual multiframelets can always be derived from a pair of compactly supported refinable vector functions through the OEP. Next we will introduce the concept of balanced moment correction filters (see Definition~\ref{def:mcf}) to characterize the structures of OEP-based dual multiframelets with high balancing orders. We will present and prove our main results Theorem~\ref{thm:mcf} of this paper in this section. In Section~\ref{sec:exmp}, we will provide an example of balanced moment correction filters to illustrate our main results.

\section{Preliminaries}\label{sec:pre}

To better explain our motivations and to get prepared for our later study, we provide all necessary backgrounds on dual framelets and the OEP in this section.

\subsection{Backgrounds on Dual Framelets}
First, let us recall some basic concepts. Throughout this paper, $\dm$ is a $d\times d$ dilation matrix, i.e., $\dm\in\Z^{d\times d}$ and its eigenvalues are all greater than one in modulus. For simplicity, let 
$$d_{\dm}:=|\det(\dm)|.$$ 
Denote $ (\dLp{2})^{r\times s}$ the linear space of all $r\times s$ matrices of square integrable functions in $\dLp{2}$. For simplicity, $(\dLp{2})^r:=(\dLp{2})^{r\times 1}$. We introduce the following notion:
$$\la f,g\ra:=\int_{\dR}f(x)\ol{g(x)}^{\tp}dx,\qquad \forall f\in\dLrs{2}{r}{s},\quad g\in\dLrs{2}{t}{s}.$$
Let $\mrphi,\mrtphi\in (\dLp{2})^r$, $\psi, \tpsi\in (\dLp{2})^s$. We say that $\{\mrphi;\psi\}$ is \emph{an $\dm$-framelet} in $\dLp{2}$ if
there exist positive constants $C_1$ and $C_2$ such that
%
$$C_1\|f\|_{\dLp{2}}^2\le
\sum_{k\in \dZ} |\la f, \mrphi(\cdot-k)\ra|^2+\sum_{j=0}^\infty \sum_{k\in \dZ}
|\la f, \psi_{\dm^j;k}\ra|^2\le C_2\|f\|_{\dLp{2}}^2, \quad  f\in \dLp{2},$$
where $\psi_{\dm^j;k}:=d_{\dm}^{j/2}\psi(\dm^j\cdot-k)$ and $|\la f, \psi_{\dm^j;k}\ra|^2:=\|\la f, \psi_{\dm^j;k}\ra\|^2_{l_2}$. $(\{\mrphi; \psi\},\{\mrtphi;\tpsi\})$ is called \emph{a dual $\dm$-framelet} in $\dLp{2}$ if both $\{\mrphi; \psi\}$ and $\{\mrtphi; \tpsi\}$ are $\dm$-framelets in $\dLp{2}$ and satisfy
\be\label{qtf:expr}
f=\sum_{k\in \dZ} \la f, \mrphi(\cdot-k)\ra \mrtphi(\cdot-k)+
\sum_{j=0}^\infty \sum_{k\in \dZ}
\la f, \psi_{\dm^j;k}\ra \tpsi_{\dm^j;k}, \qquad \forall f\in \dLp{2},
\ee
with the above series converging unconditionally in $\dLp{2}$. $(\{\mrphi; \psi\},\{\mrtphi;\tpsi\})$ is called \emph{a dual multiframelet} if the multiplicity $r>1$, and is called \emph{a scalar dual framelet} if $r=1$. Unless specified, we shall use the term framelet to refer to both.

For a dual $\dm$-framelet  $(\{\mrphi;\psi\}, \{\tilde{\mrphi}; \tilde{\psi}\})$, the sparseness of the frame expansion \eqref{qtf:expr} is closely related to the vanishing moments on the framelet generators $\psi$ and $\tpsi$. We say that $\psi$ has $m$ vanishing moments if
$$\int_{\dR}\pp(x)\psi(x)dx=0,\qquad \forall \pp\in\PL_{m-1},$$
where $\PL_{m-1}$ is the space of all $d$-variate polynomials of degree at most $m-1$. Note that $\psi$ has $m$ vanishing moments if and only if
$$\wh{\psi}(\xi)=\bo(\|\xi\|^m),\qquad \xi\to 0,$$
where $f(\xi)=g(\xi)+\bo(\|\xi\|^m)$ as $\xi\to 0$ means $\partial^\mu f(0)=\partial^\mu g(0)$ for all $\mu\in\dN_{0;m}$ with $|\mu|:=\mu_1+\cdots+\mu_d<m$.
We define $\vmo(\psi):=m$ with $m$ being the largest such integer. It is well known in approximation theory (see e.g. \cite[Proposition 5.5.2]{hanbook}) that if $\vmo(\psi)=m$ and $\vmo(\tpsi)=\tilde{m}$, then we necessarily have 
\be \label{qi:order}
\sum_{k\in \dZ} \la \pp, \mrphi(\cdot-k)\ra \mrtphi(\cdot-k):=\sum_{\ell=1}^r \sum_{k\in \dZ} \la \pp, \mrphi_\ell(\cdot-k)\ra \mrtphi_\ell(\cdot-k)=\pp,\qquad \forall\, \pp\in \PL_{m-1}.
\ee
which plays a crucial role in approximation theory and numerical analysis for the convergence rate of the associated approximation/numerical scheme. Moreover, we have
$$\ol{\wh{\mrphi}(\xi)}^\tp \wh{\mrtphi}(\xi+2\pi k)=\bo(\|\xi\|^{m}), \qquad\ol{\wh{\mrphi}(\xi+2\pi k)}^\tp \wh{\mrtphi}(\xi)=\bo(\|\xi\|^{\tilde{m}}),\quad k\in\dZ\bs\{0\},$$
and
$$\ol{\wh{\mrphi}(\xi)}^\tp  \wh{\mrtphi}(\xi)=1+\bo(\|\xi\|^{\tilde{m}+m}),$$
as $\xi\to 0$.

A popular method called the \emph{oblique extension principle (OEP)} has been introduced in the literature, which allows us to construct dual framelets with all generators having sufficiently high vanishing moments from \emph{refinable vector functions} \cite{cpss13,cpss15,cs08,ch01,cj00,dhacha,fjs16,hjsz18,js15,kps16,lj06,sz16}. Denote $\dlrs{0}{r}{s}$ the linear space of all $r\times s$ matrix-valued sequences $u=\{u(k)\}_{k\in \dZ}:\dZ\to \C^{r\times s}$ with finitely many non-zero terms. Any element $u\in\dlrs{0}{r}{s}$ is said to be a \emph{finitely supported (matrix-valued) filter/mask}. For $\phi\in (\dLp{2})^r$, we say that $\phi$ is \emph{an $\dm$-refinable vector function} with \emph{a refinement filter/mask} $a\in \dlrs{0}{r}{r}$ if the following \emph{refinement equation} is satisfied:
%
\be\label{ref}\phi(x)=d_{\dm}\sum_{k\in\dZ}a(k)\phi(\dm x-k),\qquad   x\in\dR.\ee
If $r=1$, then we simply say that $\phi$ is an $\dm$-refinable (scalar) function. For $u\in\dlrs{0}{r}{s}$, define its Fourier series via $\wh{u}(\xi):=\sum_{k\in\dZ}u(k)e^{-ik\cdot \xi}$ for $\xi\in\dR$. The Fourier transform is defined via $\wh{f}(\xi):=\int_{\dR} f(x) e^{-ix\cdot\xi} dx$ for $\xi\in \dR$ for all $f\in \dLp{1}$, and can be naturally extended to $\dLp{2}$ functions and tempered distributions. The refinement equation \eqref{ref} is equivalent to
\be\label{ref:f}\wh{\phi}(\dm^{\tp} \xi)=\wh{a}(\xi)\wh{\phi}(\xi),\qquad \xi\in \dR,\ee
where $\wh{\phi}$ is the $r\times 1$ vector obtained by taking entry-wise Fourier transform on $\phi$. Most known framelets are constructed from refinable vector functions via the OEP, and we refer to them as OEP-based framelets. There are several versions of the OEP which have been introduced in the literature (see \cite{chs02,dhrs03,hanbook,hl20pp}). Here we recall the following version of the OEP for compactly supported multivariate multiframelets:
\begin{theorem}[\emph{Oblique extension principle (OEP)}]\label{thm:df}
	Let $\dm$ be a $d\times d$ dilation matrix. Let $\theta,\tilde{\theta},a,\tilde{a}\in \dlrs{0}{r}{r}$ and $\phi,\tilde{\phi}\in (\dLp{2})^r$ be compactly supported $\dm$-refinable vector functions with refinement filters $a$ and $\tilde{a}$, respectively.
	For matrix-valued filters $b,\tilde{b}\in \dlrs{0}{s}{r}$,
	define
	\be\label{oep:mr}\wh{\mrphi}(\xi):=\wh{\theta}(\xi)\wh{\phi}(\xi),
	\quad \wh{\psi}(\xi):=\htb(\dm^{-\tp}\xi)\wh{\phi}(\dm^{-\tp}\xi),\ee
	\be\label{oep:mrt}\wh{\tilde{\mrphi}}(\xi):=\wh{\tilde{\theta}}(\xi)\wh{\tilde{\phi}}(\xi),
	\quad \wh{\tilde{\psi}}(\xi):=\wh{\tilde{b}}(\dm^{-\tp}\xi)\wh{\tilde{\phi}}(\dm^{-\tp}\xi).
	\ee
	Then $(\{\mrphi;\psi\}, \{\tilde{\mrphi}; \tilde{\psi}\})$ is a dual $\dm$-framelet in $\dLp{2}$ if the following conditions are satisfied:
	\begin{enumerate}
		\item[(1)] $\ol{\wh{\phi}(0)}^\tp \wh{\Theta}(0)\wh{\tilde{\phi}}(0)=1$ with $\wh{\Theta}(\xi):=\ol{\wh{\theta}(\xi)}^\tp \wh{\tilde{\theta}}(\xi)$;
		
		\item[(2)]  $\wh{\psi}(0)=\wh{\tilde{\psi}}(0)=0$.
		
		\item[(3)] $(\{a;b\},\{\tilde{a};\tilde{b}\})_{\Theta}$ forms an OEP-based dual $\dm$-framelet filter bank, i.e.,
		\be \label{dffb}
		\ol{\wh{{a}}(\xi)}^\tp \wh{\Theta}(\dm^{\tp} \xi){\wh{\tilde{a}}(\xi+2\pi \omega)}+\ol{\wh{{b}}(\xi)}^\tp {\wh{\tilde{b}}(\xi+2\pi \omega)}=\td(\omega)\wh{\Theta}(\xi),\ee
		for all $\xi\in\dR$ and $\omega\in\Omega_{\dm}$, where
		\be \label{delta:seq}
		\td(0):=1 \quad \mbox{and}\quad
		\td(x):=0,\qquad \forall\, x\ne 0
		\ee
		and $\Omega_\dm$ is a particular choice of the representatives of cosets in $[\dm^{-\tp}\dZ]/\dZ$ given by
		\be\label{omega:dm}
		\Omega_{\dm}:=\{\om{1},\dots,\om{d_{\dm}}\}:=(\dm^{-\tp}\dZ)\cap [0,1)^d\quad \mbox{with}\quad \om{1}:=0.
		\ee
	\end{enumerate}
\end{theorem}
The key step to construct an OEP-based dual framelet is to obtain filters $\theta,\tilde{\theta}\in\dlrs{0}{r}{r}$ and $b,\tilde{b}\in\dlrs{0}{s}{r}$ such that $(\{a;b\},\{\tilde{a};\tilde{b}\})_{\Theta}$ is a dual framelet filter bank which satisfies \eqref{dffb}. For any $u\in\dlrs{0}{s}{r}$, define 
\be \label{Pb}
P_{u;\dm}(\xi):=[\wh{u}(\xi+2\pi\om{1}), \ldots, \wh{u}(\xi+2\pi \om{d_{\dm}})],\qquad\xi\in\dR,
\ee
which is an $s\times (rd_{\dm})$ matrix of $2\pi\dZ$-periodic $d$-variate trigonometric polynomials. It is obvious that \eqref{dffb} is equivalent to
\be \label{spectral}
\ol{P_{b;\dm}(\xi)}^\tp P_{\tilde{b};\dm}(\xi)=
\cM_{a,\tilde{a},\Theta}(\xi),
\ee
where 
\be\label{m:a:ta}
\begin{aligned}\cM_{a,\tilde{a},\Theta}(\xi):=&\DG\left(\wh{\Theta}(\xi+2\pi\om{1}),\ldots,
	\wh{\Theta}(\xi+2\pi \om{d_{\dm}})\right)\\
	&-\ol{P_{a;\dm}(\xi)}^\tp \wh{\Theta}(\dm^\tp\xi) P_{\tilde{a};\dm}(\xi).\end{aligned}\ee

For an OEP-based dual $\dm$-framelet $(\{\mrphi;\psi\}, \{\tilde{\mrphi}; \tilde{\psi}\})$, The orders of vanishing moments of $\psi$ and $\tpsi$ are closely related to the sum rules of the filters $a$ and $\tilde{a}$ associated to $\phi$ and $\tphi$. We say that a filter $a\in (\dlp{0})^{r\times r}$ has \emph{order $m$ sum rules with respect to $\dm$} with a matching filter $\vgu\in \dlrs{0}{1}{r}$ if $\wh{\vgu}(0)\ne 0$ and
\be \label{sr}
\wh{\vgu}(\dm^{\tp}\xi)\wh{a}(\xi+2\pi\omega)=\td(\omega)\wh{\vgu}(\xi)+
\bo(\|\xi\|^m),\quad \xi\to 0,\quad\forall\, \omega\in\Omega_{\dm}.
\ee
In particular, we define 
$$\sr(a,\dm):=\sup\{m\in\N_0: \text{\eqref{sr} holds for some }v\in\dlrs{0}{1}{r}\}.$$
It can be easily deduced from \eqref{dffb} that $\vmo(\psi)\le \sr(\tilde{a},\dm)$ and $\vmo(\tpsi)\le \sr(a,\dm)$ always hold no matter how we choose $\theta$ and $\tilde{\theta}$. Therefore, we are curious about whether or not one can construct filters $\theta,\tilde{\theta}\in\dlrs{0}{r}{r}$ in a way such that the matrix $\cM_{a,\tilde{a},\Theta}$ admits a factorization as in \eqref{spectral} for some $b,\tilde{b}\in\dlrs{0}{s}{r}$ such that $\wh{b}(\xi)\wh{\phi}(\xi)=\bo(\|\xi\|^{\sr(\tilde{a},\dm)})$ and $\wh{\tilde{b}}(\xi)\wh{\tphi}(\xi)=\bo(\|\xi\|^{\sr(a,\dm)})$ as $\xi\to 0$.
With OEP, a lot of compactly supported scalar dual framelets with the highest possible vanishing moments have been constructed in the literature, to mention only a few, see \cite{cs10,ch00,chs02,dh04,dhrs03,dhacha,dh18pp,han97,han09,hanbook,hm03,hm05,jqt01,js15,mothesis,rs97,sel01} and many references therein.

\subsection{The Perfect Reconstruction Property and the Balancing Property of a DmFT}

We now review some basic concepts and properties of a DmFT. By $(\dsq)^{s\times r}$ we denote the linear space of all sequences $v: \dZ \rightarrow \C^{s\times r}$. We call every element $v\in(\dsq)^{s\times r}$ \emph{a matrix-valued filter}. For a filter $a\in \dlrs{0}{r}{r}$, we define the filter $a^\star$ via $\wh{a^\star}(\xi):=\ol{\wh{a}(\xi)}^\tp$, or equivalently,
$a^\star(k):=\ol{a(-k)}^\tp$ for all $k\in \dZ$.
We define the \emph{convolution} of two filters via
\[
[v*u](n):=\sum_{k\in \Z} v(k) u(n-k),\quad n\in \dZ,\quad v\in(\dsq)^{s\times r},\quad u\in\dlrs{0}{r}{t}.
\]
Let $\dm$ be a $d\times d$ dilation matrix, define the \emph{upsampling operator} $\uparrow \dm: (\dsq)^{s\times r}\to(\dsq)^{s\times r}$ as
$$[v\uparrow \dm](k):=\begin{cases}v(\dm^{-1} k), &\text{if }k\in \dZ\cap[\dm^{-1}\dZ],\\
0, &\text{elsewhere},\end{cases},\qquad  \forall k\in\dZ,\quad v\in(\dsq)^{s\times r}.$$

We introduce the following operators acting on matrix-valued sequence spaces:

\begin{itemize} 
	
	\item For $u\in\dlrs{0}{r}{t}$, the \emph{subdivision operator} $\sd_{u,\dm}$ is defined via
	$$\sd_{u,\dm} v=|\det(\dm)|^{\frac{1}{2}}[v\uparrow\dm]*u=|\det(\dm)|^{\frac{1}{2}}\sum_{k\in\dZ}v(k)u(\cdot-\dm k),$$
	for all $v\in(\dsq)^{s\times r}$.
	
	\item For $u\in\dlrs{0}{t}{r}$, the \emph{transition operator} $\tz_{u,\dm}$ is defined via
	$$\tz_{u,\dm} v=|\det(\dm)|^{\frac{1}{2}}[v*u^\star]\downarrow\dm =|\det(\dm)|^{\frac{1}{2}}\sum_{k\in\dZ}v(k)\ol{u(k-\dm\cdot)}^{\tp},$$
	for all $v\in(\dsq)^{s\times r}$.
\end{itemize}

Let $\theta,\tilde{\theta},a,\tilde{a}\in\dlrs{0}{r}{r}$ and $b,\tilde{b}\in\dlrs{0}{s}{r}$ be finitely supported filters. For any $J\in\N$ and any input data $v_0\in(\dsq)^{1\times r}$, the \emph{$J$-level discrete multiframelet transform (DmFT)} employing the filter bank $(\{a;b\},\{\tilde{a};\tilde{b}\})_{\Theta}$ where $\Theta:=\theta^{\star}*\tilde{\theta}$ is implemented as follows:
\begin{enumerate}
	\item[(S1)] \emph{Decomposition/Analysis}: Recursively compute $v_j,w_j$ for $j=1,\dots,s$ via
	\be\label{dft:anal}v_j:=\tz_{a,\dm}v_{j-1},\qquad w_j:=\tz_{b,\dm}v_{j-1}.\ee
	
	\item[(S2)] \emph{Reconstruction/Synthesis}: Define $\tilde{v}_J:=v_J*\Theta$. Recursively compute $\tilde{v}_{j-1}$ for $j=J,\dots,1$ via
	\be\label{dft:syn}\tilde{v}_{j-1}:=\sd_{\tilde{a},\dm}\mrv_j+\sd_{\tilde{b},\dm}w_j.\ee
	
	\item[(S3)] \emph{Deconvolution}: Recover $\breve{v}_0$ from $\tilde{v}_0$ through $\breve{v}_0*\Theta=\tilde{v}_0.$
\end{enumerate}

We call $\{a;b\}$ the \emph{analysis filter bank} and $\{\tilde{a};\tilde{b}\}$ the \emph{synthesis filter bank}. If any input data $v\in(\dsq)^{1\times r}$ can be exactly retrieved from the above transform, then we say that the $J$-level DmFT has the \emph{perfect reconstruction (PR) property}. By \cite[Theorem 2.3]{hl20pp}, a $J$-level DmFT has the PR property for all $J\in\N$ if and only if $(\{a;b\},\{\tilde{a},\tilde{b}\})_{\Theta}$ is a dual $\dm$-framelet filter bank and $\Theta$ is \emph{strongly invertible}. Here we say that $\Theta\in\dlrs{0}{r}{r}$ is strongly invertible if if there exists $\Theta^{-1}\in\dlrs{0}{r}{r}$ such that $\wh{\Theta^{-1}}=\wh{\Theta}^{-1}$, or equivalently $\det(\wh{\Theta}(\xi))=ce^{-ik\cdot\xi}$ for some $c\in\C\bs\{0\}$ and $k\in\Z^d$. Note that if $\Theta$ is strongly invertible, then so are both $\theta$ and $\tilde{\theta}$ as $\wh{\Theta}=\ol{\wh{\theta}}^\tp\wh{\tilde{\theta}}$.

In many applications, the original data is scalar-valued, that is, an input data $v\in\dsq$. Thus to implement a multi-level DmFT, we need to first vectorize the input data. Let $\dn$ be a $d\times d$ integer matrix with $|\det(\dn)|=r$, and let $\Gamma_{\dn}$ be a particular choice of the representatives of the cosets in $\dZ/[\dn\dZ]$ given by 
\be\label{Gamma:N}\Gamma_{\dn}:=\{\mrgam_1,\dots,\mrgam_r\}=:[\dn[0,1)^d]\cap\dZ,\quad\mbox{with}\quad \mrgam_1:=0.\ee
We define the \emph{standard vectorization operator with respect to $\dn$} via
\be\label{vec:con}E_{\dn}v:=(v(\dn\cdot+\mrgam_1),\dots v(\dn\cdot+\mrgam_r)),\qquad\forall v\in\dsq.\ee
Clearly $E_{\dn}$ is a bijection between $\dsq$ and $(\dsq)^{1\times r}$. The sparsity of a multi-level DmFT employing a dual $\dm$-framelet filter bank $(\{a;b\},\{\tilde{a};\tilde{b}\})_{\Theta}$ ($\Theta:=\theta^\star*\tilde{\theta}$) is measured by the \emph{$E_{\dn}$-balancing order} of the analysis filter bank $\{a;b\}$, denoted by $\bpo(\{a;b\},\dm,\dn):=m$ where $m$ is the largest integer such that the following two conditions hold:
\begin{enumerate}
	\item[(i)] $\tz_{a,\dm}$ is invariant on $E_{\dn}(\PL_{m-1}|_{\dZ})$, i.e., \be \label{bp:lowpass}
	\tz_{a,\dm} E_{\dn}(\PL_{m-1}|_{\dZ})\subseteq E_{\dn}(\PL_{m-1}|_{\dZ}).\ee
	\item[(ii)] The filter $b$ has $m$ \emph{$E_{\dn}$-balancing vanishing moments},  i.e.,
	\be \label{bp:highpass}
	\tz_{b,\dm} E_{\dn} (\pp)=0, \qquad \forall\pp\in \PL_{m-1}|_{\dZ}.
	\ee
	
\end{enumerate}
If items (i) and (ii) are satisfied, note that the framelet coefficient $w_j=\tz_{b,\dm}\tz_{a,\dm}^{j-1}E_{\dn}(\pp)=0$ for all $\pp\in\PL_{m-1}|_{\dZ}$ and $j=1,\dots,J$. This preserves sparsity at all levels of the multi-level DmFT. A complete characterization of the balancing order of a filter bank is given by the following result.

\begin{theorem}\cite[Proposition~3.1, Theorem 4.1]{han10}\label{thm:bp}
	Let $\dm$ be a $d \times d$ dilation matrix and $r\ge 2$ be a positive integer.
	Let $a\in \dlrs{0}{r}{r}$ and $b\in \dlrs{0}{s}{r}$ for some $s\in\N$. Let $\dn$ be a $d\times d$ integer matrix with $|\det(\dn)|=r$ and $E_{\dn}$ in \eqref{vec:con}.
	Define 
	\be \label{vgu:special}
	\wh{\Vgu_{\dn}}(\xi):=\left(e^{i\dn^{-1}\ka{1}\cdot\xi},\ldots, e^{i \dn^{-1}\ka{r}\cdot\xi}\right),\qquad\xi\in\dR.
	\ee
	Then the following statements hold:
	\begin{enumerate}
		
		\item[(1)] The filter $b$ has order $m$ $E_{\dn}$-balancing vanishing moments satisfying \eqref{bp:highpass} if and only if
		\be \label{cond:bvmo} \wh{\Vgu_{\dn}}(\xi)\ol{\wh{b}(\xi)}^\tp
		=\bo(\|\xi\|^m),\qquad \xi \to 0.
		\ee
		\item[(2)] The filter bank $\{a;b\}$ has $m$ $E_{\dn}$-balancing order if and only if \eqref{cond:bvmo} holds and
		%
		$$\wh{\Vgu_{\dn}}(\xi)\ol{\wh{a}(\xi)}^\tp
		=	 \wh{c}(\xi)\wh{\Vgu_{\dn}}(\dm^{\tp}\xi)
		+\bo(\|\xi\|^m),\quad \xi\to 0,$$
		for some $c\in \dlp{0}$ with $\wh{c}(0)\ne 0$.
		
	\end{enumerate}
\end{theorem} 

Suppose that $\phi,\tilde{\phi}\in (\dLp{2})^r$ are compactly supported $\dm$-refinable vector functions in $\dLp{2}$ satisfying $\wh{\phi}(\dm^{\tp}\xi)=\wh{a}(\xi)\wh{\phi}(\xi)$
and
$\wh{\tilde{\phi}}(\dm^{\tp}\xi)=\wh{\tilde{a}}(\xi)\wh{\tilde{\phi}}(\xi)$.
Define $\mrphi,\psi,\mrtphi,\tpsi$ as in \er{oep:mr} and \er{oep:mrt}.
If $\ol{\wh{\phi}(0)}^\tp \wh{\Theta}(0)\wh{\tilde{\phi}}(0)=1$ and $\wh{\psi}(0)=\wh{\tilde{\psi}}(0)=0$,
then Theorem~\ref{thm:df} tells us that $(\{\mrphi;\psi\},\{\tilde{\mrphi};\tilde{\psi}\})$ is a dual $\dm$-framelet in $\dLp{2}$. With $m:=\sr(\tilde{a},\dm)$, we observe that $ \vmo(\psi)\le m$, $\bvmo(b,\dm,\dn)\le m$ and $\bpo(\{a;b\},\dm,\dn)\le \bvmo(b,\dm,\dn)$.
If $\bpo(\{a;b\},\dm,\dn)=\bvmo(b,\dm,\dn)
=\vmo(\psi)=m$, then we say that the DmFT (or the dual multiframelet $(\{\mrphi;\psi\},\{\tilde{\mrphi};\tilde{\psi}\})$) is \emph{$E_{\dn}$-balanced}.
For $r>1$,  $\bpo(\{a,b\},\dm,\dn)<\vmo(\psi)$ often happens. Hence, having high vanishing moments on framelet generators does not guarantee the balancing property and thus significantly reduces the sparsity of the associated discrete multiframelet transform. How to overcome this shortcoming has been extensively studied in the setting of functions in \cite{cj00,lv98,sel00} and in the setting of DmFT in \cite{han09,han10,hanbook}.\\

\subsection{The normal form of a matrix-valued filter}

In this subsection, we briefly review the results of a recently developed normal form of the matrix-valued filter. The matrix-valued filter normal form greatly reduces the difficulty in studying multiframelets and multiwavelets, in a way such that we can mimic the techniques we have for studying scalar framelets and wavelets. Considerable works on this topic have been done. We refer the readers to  \cite{han03,han09,han10,hanbook,hl20pp,hm03} for detailed discussion. The most recent advance on this topic is \cite{hl20pp}, which not only generalizes all previously existing works under much weaker conditions but also provides a strengthened normal form of a matrix-valued filter which greatly benefits our study on balanced multivariate multiframelets.

We first recall the following lemma which is known as \cite[Lemma 2.2]{han10}. This result links different vectors of functions which are smooth at the origin by strongly invertible filters.

\begin{lemma}\label{lem:U}[\cite[Lemma 2.2]{han10}]
	Let $\wh{v}=(\wh{v_1},\ldots,\wh{v_r})$ and $\wh{u}=(\wh{u_1},\ldots,\wh{u_r})$ be $1\times r$ vectors of functions which are infinitely differentiable at $0$ with $\wh{v}(0)\ne 0$ and $\wh{u}(0)\ne 0$. If $r\ge 2$, then for any positive integer $n\in \N$, there exists a strongly invertible $U\in\dlrs{0}{r}{r}$ such that
	$\wh{u}(\xi)=\wh{v}(\xi)\wh{U}(\xi)+\bo(\|\xi\|^n)$ as $\xi\to 0$.
\end{lemma}

One of the most important results on the normal form of a matrix-valued filter is the following result which has been developed recently and is a part of \cite[Theorem 3.3]{hl20pp}.

\begin{theorem}\label{thm:normalform:gen}Let $\wh{v},\wh{\mrv}$ be $1\times r$ vectors and $\wh{\phi},\wh{\mrphi}$ be $r\times 1$ vectors of functions which are infinitely differentiable at $0$. Suppose 
	$$\wh{v}(\xi)\wh{\phi}(\xi)=1+\bo(\|\xi\|^m)\quad\mbox{and}\quad \wh{\mrv}(\xi)\wh{\mrphi}(\xi)=1+\bo(\|\xi\|^m),\quad\xi\to 0.$$
	If $r\ge 2$, then for each $n\in\N$, there exists a strongly invertible filter $U\in\dlrs{0}{r}{r}$ such that
	$$\wh{v}(\xi)\wh{U}(\xi)^{-1}=\wh{\mrv}(\xi)+\bo(\|\xi\|^m)\quad\mbox{and}\quad \wh{U}(\xi)\wh{\phi}(\xi)=\wh{\mrphi}(\xi)+\bo(\|\xi\|^n),\quad\xi\to 0.$$
	
\end{theorem}

A special case of Theorem~\ref{thm:normalform:gen} is the following result (\cite[Theorem 1.2]{hl20pp}, cf. \cite[Theorem 5.1]{han10}).

\begin{theorem}\label{thm:normalform}Let $\dm$ be a $d\times d$ dilation matrix, and let $m\in\N$ and $ r\ge 2$ be integers.
	Let $\phi$ be an $r\times 1$ vector of
	compactly supported distributions satisfying $\wh{\phi}(\dm^{\tp} \xi)=\wh{a}(\xi)\wh{\phi}(\xi)$ with $\wh{\phi}(0)\ne 0$ for some $a\in\dlrs{0}{r}{r}$.
	Suppose the filter $a$ has order $m$ sum rules with respect to $\dm$ satisfying \eqref{sr} with a matching filter $\vgu\in \dlrs{0}{1}{r}$ such that $\wh{\vgu}(0)\wh{\phi}(0)=1$.
	Then for any positive integer $n\in \N$, there exists a strongly invertible filter $U\in\dlrs{0}{r}{r}$ such that the following statements hold:
	\begin{enumerate}
		\item[(1)] Define
		$\wh{\mathring{\vgu}}:=(\wh{\mathring{\vgu}_1},\ldots,	 \wh{\mathring{\vgu}_r}):=\wh{\vgu}\wh{U}^{-1}$ and
		$\wh{\mathring{\phi}}:=(\wh{\mathring{\phi}_1},\ldots,		 \wh{\mathring{\phi}_r})^\tp:=\wh{U}\wh{\phi}$.
		We have
		\be \label{normalform:phi}		 \wh{\mathring{\phi}_1}(\xi)=1+\bo(\|\xi\|^n)
		\quad \mbox{and}\quad \wh{\mathring{\phi}_\ell}(\xi)=\bo(\|\xi\|^n),\quad \xi\to 0,\quad \ell=2,\ldots,r,
		\ee		 \be\label{normalform:vgu}\wh{\mathring{\vgu}_1}(\xi)=1+\bo(\|\xi\|^m)
		\quad \mbox{and}\quad
		\wh{\mathring{\vgu}_\ell}(\xi)=\bo(\|\xi\|^m),\quad \xi\to 0,\quad \ell=2,\ldots,r.\ee	
		
		\item[(2)] Define $\mra\in\dlrs{0}{r}{r}$ via  $\wh{\mathring{a}}:=\wh{U}(\dm^{\tp}\cdot) \wh{a}\wh{U}^{-1}$. Then $\wh{\mrphi}(\dm^{\tp}\cdot)=\wh{\mra}\wh{\mrphi}$ and the new filter $\mra$ has order $m$ sum rules with respect to $\dm$ with the matching filter $\mathring{\vgu}\in (\dlp{0})^{1\times r}$.

	\end{enumerate}
	
\end{theorem}

Let $\mra\in\dlrs{0}{r}{r}$ be a refinement mask associated to an $\dm$-refinable vector function $\mrphi$ satisfying \er{normalform:phi}, and suppose that $\mra$ has $m$ sum rules with a matching filter $\mrvgu\in\dlrs{0}{1}{r}$ satisfying \er{normalform:vgu}. It is not hard to observe that $\mra$ has the following structure:
\be \label{normalform}
\wh{\mra}(\xi)=\left[
\begin{matrix}\wh{\mra_{1,1}}(\xi) & \wh{\mra_{1,2}}(\xi)\\
	\wh{\mra_{2,1}}(\xi) & \wh{\mra_{2,2}}(\xi)\end{matrix}\right],
\ee
where $\wh{\mra_{1,1}},\wh{\mra_{1,2}}, \wh{\mra_{2,1}}$ and $\wh{\mra_{2,2}}$ are $1\times 1$, $1\times (r-1)$, $(r-1)\times 1$ and $(r-1)\times (r-1)$ matrices of $2\pi\dZ$-periodic trigonometric polynomials such that
\begin{align}	 &\wh{\mra_{1,1}}(\xi)=1+\bo(\|\xi\|^n),\quad \wh{\mra_{1,1}}(\xi+2\pi\omega)=\bo(\|\xi\|^m),\quad \xi\to 0,\quad\forall \omega\in\Omega_{\dm}\setminus\{0\},\label{mra:11}\\
&\wh{\mra_{1,2}}(\xi+2\pi\omega)=\bo(\|\xi\|^m),\quad \xi\to 0,\quad\forall \omega\in\Omega_{\dm},\label{mra:12}\\ &\wh{\mra_{2,1}}(\xi)=\bo(\|\xi\|^n),\quad \xi\to 0.\label{mra:21}
\end{align}
Any filter $\mra$ satisfying \er{normalform}, \er{mra:11}, \er{mra:12} and \er{mra:21} is said to take the \emph{ideal $(m,n)$-normal form}.\\

If $d=1$, then the three moment conditions \er{mra:11}, \er{mra:12} and \er{mra:21} further yield
$$\wh{\mra_{1,1}}(\xi)=(1+e^{-i\xi}+\dots+e^{-i(|\dm|-1)\xi})^{m}P_{1,1}(\xi)=1+\bo(|\xi|^n),\qquad\xi\to 0,$$
$$\wh{\mra_{1,2}}(\xi)=(1-e^{-i|\dm|\xi})^mP_{1,2}(\xi),\qquad  \wh{\mra_{2,1}}(\xi)=(1-e^{-i\xi})^nP_{2,1}(\xi),$$
where $P_{1,1}, P_{1,2}$ and $P_{2,1}$ are some $1\times 1, 1\times (r-1)$ and $(r-1)\times 1$ matrices of $2\pi$-periodic trigonometric polynomials. Recall that a $2\pi$-periodic trigonometric polynomial $\wh{u}$ satisfies $\wh{u}(\xi)=\bo(|\xi|^m)$ as $\xi\to 0$ if and only if $(1-e^{-i\xi})^m$ divides $\wh{u}$. This is the crucial property to construct univariate dual framelets with high vanishing moments. Unfortunately for $d\ge 2$, there are no corresponding factors for $(1+e^{-i\xi}+\dots+e^{-i(|\dm|-1)\xi})^m$ and $(1-e^{-i\xi})^m$. This means the factorization technique that we have to construct dual framelets with high vanishing moments for the case $d=1$ is no longer available, which illustrates that the investigation is more difficult for $d>1$.

\section{The Existence and a Structural Characterization of Balanced Dual Multiframelets}\label{sec:cha}
In this section, we shall establish our main results on the structural characterization of balanced dual multiframelets. First, we will prove in Theorem~\ref{thm:dfrt} that a balanced dual multiframelet can always be constructed from an arbitrary pair of compactly supported refinable vector functions. Motivated by Theorem~\ref{thm:dfrt}, the concept of balanced moment correction filters will be introduced, which is exactly what we are going to use to study the structure of balanced dual multiframelets, and the main result of this part is summarized as Theorem~\ref{thm:mcf}.

\subsection{The Existence of OEP-based Balanced Dual Multiframelets in $\R^d$}. The first question to be answered is whether it is always possible to construct a balanced dual multiframelet from an arbitrary pair of refinable vector functions, which leads us to the following main theorem.

\begin{theorem}\label{thm:dfrt} Let $\dm$ be a $d\times d$ dilation matrix and $r\ge 2$ be an integer. Let $\phi,\tilde{\phi}\in(\dLp{2})^r$ be compactly supported $\dm$ refinable vector functions associated with refinement masks $a,\tilde{a}\in\dlrs{0}{r}{r}$. Suppose that $\sr(a,\dm)=\tilde{m}$ and $\sr(\tilde{a},\dm)=m$ with matching filters $\vgu,\tvgu\in\dlrs{0}{1}{r}$ respectively such that $\wh{\vgu}(0)\wh{\phi}(0)\ne 0$ and $\wh{\tvgu}(0)\wh{\tphi}(0)\ne 0$. Let $\dn$ be a $d\times d$ integer matrix with $|\det(\dn)|=r$. Then there exist $\theta,\tilde{\theta}\in\dlrs{0}{r}{r}$ and $b,\tilde{b}\in\dlrs{0}{s}{r}$ for some $s\in\N$ such that
	
	\begin{enumerate}
		\item[(1)] $\theta$ and $\tilde{\theta}$ are both strongly invertible.
		
		\item[(2)] Define finitely supported filters $\mra,\mrb,\mrta,\mrtb$ via 
		
		\be\label{mr:filter}\wh{\mra}(\xi)=\wh{\theta}(\dm^\tp \xi)\wh{a}(\xi)\wh{\theta}(\xi)^{-1},\quad \wh{\mrb}(\xi)=\wh{b}(\xi)\wh{\theta}(\xi)^{-1},\ee
		\be\label{mrt:filter}\wh{\mrta}(\xi)=\httheta(\dm^\tp \xi)\hta(\xi)\httheta(\xi)^{-1},\quad \wh{\mrtb}(\xi)=\htb(\xi)\httheta(\xi)^{-1}.\ee Then $(\{\mra;\mrb\},\{\mrta;\mrtb\})_{\td I_r}$ is an OEP-based dual $\dm$-framelet filter bank. Moreover, the discrete framelet transform employing the filter bank $(\{\mra;\mrb\},\{\mrta;\mrtb\})_{\td I_r}$ is order $m$ $E_{\dn}$-balanced, i.e.,  $\bpo(\{\mra;\mrb\},\dm,\dn)=\sr(\mrta,\dm)=m$.

		\item[(3)] $(\{\mrphi;\psi\},\{\mrtphi;\tpsi\})$ is a compactly supported dual $\dm$-framelet in $\dLp{2}$ with $\vmo(\psi)=\tilde{m}$ and $\vmo(\tpsi)=m$, where $\mrphi,\psi,\mrtphi,\tpsi$ are vector-valued functions defined as in \er{oep:mr} and \er{oep:mrt}.
	\end{enumerate}
	
\end{theorem}

\begin{definition}\label{def:mcf}The strongly invertible filters $\theta,\tilde{\theta}\in \dlrs{0}{r}{r}$ in Theorem~\ref{thm:dfrt} are called a pair of \emph{$E_{\dn}$-balanced moment correction filters} for the refinement masks $a,\tilde{a}\in\dlrs{0}{r}{r}$.
\end{definition}

To prove Theorem~\ref{thm:dfrt}, we first need to introduce several notations. For any $k\in\dZ$, the \emph{backward difference operator} $\nabla_k$ is defined via
$$\nabla_ku(n):=u(n)-u(n-k),\qquad \forall n\in\dZ,\quad u\in (\dsq)^{t\times r}.$$
For any multi-index $\alpha:=(\alpha_1,\dots,\alpha_d)\in\dN_0$, we define
$$\nabla^\alpha:=\nabla^{\alpha_1}_{e_1}\nabla^{\alpha_2}_{e_2}\dots\nabla^{\alpha_d}_{e_d},$$
where $\{e_1,\dots,e_d\}$ is the standard basis for $\dR$. Observe that
$$\wh{\nabla^{\alpha}u}(\xi)=
\wh{\nabla^{\alpha}\td}(\xi)\wh{u}(\xi)=
(1-e^{-i\xi_1})^{\alpha_1}(1-e^{-i\xi_2})^{\alpha_2}
\cdots(1-e^{-i\xi_d})^{\alpha_d}\wh{u}(\xi),$$
for all $\xi=(\xi_1,\dots,\xi_d)^{\tp}\in\dR$ and $u\in\dlrs{0}{t}{r}$. 

For $d=1$, recall that a $2\pi$-periodic trignometric polynomial $\wh{c}$ satisfies $\wh{c}(\xi)=\bo(|\xi|^m)$ as $\xi\to 0$ if and only if $\wh{c}$ is divisible by $(1-e^{-i\xi})^m$. Though such a factorization is not available when $d>1$ and there is no factor which plays the role of $(1-e^{-i\xi})^m$ as in the univariate case, the following lemma tells us exactly how one can characterize the moments at zero of a multivariate trigonometric polynomial.

\begin{lemma} (\cite[Lemma 5]{dhacha})\label{diff:vm}Let $m\in\N$ and $\wh{c}$ be a $2\pi\dZ$-periodic $d$-variate trigonometric polynomial. Then $\wh{c}(\xi)=\bo(\|\xi\|^m)$ as $\xi\to 0$ if and only if
	$$\wh{c}(\xi)=\sum_{\alpha\in\dN_{0;m}}\wh{\nabla^\alpha\td}(\xi)\wh{c_{\alpha}}(\xi)$$
	for some $c_\alpha\in\dlp{0}$ for all $\alpha\in\dN_{0;m}$, where
	\be\label{n:0:m}\dN_{0;m}:=\{\alpha\in\dN_0:|\alpha|:=\alpha_1+\dots+\alpha_d=m\}.\ee
\end{lemma}

Next, we introduce the notion of the so-called \emph{coset sequences}. Let $\dm$ be an \emph{invertible integer matrix} and let $\gamma\in\dZ$. For any matrix-valued sequence $u\in(\dsq)^{t\times r}$, we define the \emph{$\gamma$-coset sequence} of $u$ with respect to $\dm$ via
$$u^{[\gamma]}(k):=u^{[\gamma;\dm]}(k)=u(\gamma+\dm k),\quad k\in\dZ.$$
For $u\in \dlrs{0}{t}{r}$, it is easy to see that
\be\label{f:coset}\wh{u}(\xi)=\sum_{\gamma\in \Gamma_\dm}
\wh{u^{[\gamma;\dm]}}(\dm^{\tp}\xi)e^{-i\gamma\cdot\xi},\qquad\xi\in\dR,\ee
where $\Gamma_\dm$ is a complete set of canonical representatives of the quotient group $\dZ/[\dm\dZ]$, with
\be\label{enum:gamma}
\Gamma_{\dm}:=\{\ga{1},\dots,\ga{d_{\dm}}\}=:(\dm[0,1)^d)\cap\dZ
\quad \mbox{with}\quad \ga{1}:=0.
\ee
Define $\Omega_{\dm}$ via \eqref{omega:dm}. For any filter $u\in\dlrs{0}{t}{r}$ and $\omega\in\Omega_{\dm}$, we introduce the following matrices of trigonometric polynomials associated with $u$ and $\omega$:

\begin{itemize}
	\item Define the $(td_{\dm} )\times (rd_{\dm})$ block matrix $D_{u,\omega;\dm}(\xi)$, whose $(l,k)$-th $t\times r$ blocks are given by
	\be\label{Duni}(D_{u,\omega;\dm}(\xi))_{l,k}:=
	\begin{cases}\wh{u}(\xi+2\pi\omega), &\text{if }\om{l}+\omega-\om{k}\in\dZ\\
		0, &\text{otherwise}.\end{cases}
	\ee
	
	\item Define the $(td_{\dm} )\times (rd_{\dm})$ block matrix $E_{u,\omega;\dm}(\xi)$, whose $(l,k)$-th $r\times r$ blocks are given by
	\be\label{Euni}(E_{u,\omega;\dm}(\xi))_{l,k}:=
	\wh{u^{[\ga{k}-\ga{l};\dm]}}
	(\xi)e^{-i\ga{k}\cdot(2\pi\omega)}.
	\ee
	
	\item  Define the $t\times (r d_{\dm})$ matrix $Q_{u;\dm}(\xi)$ via
	\be \label{Qu}
	Q_{u;\dm}(\xi):=\big[\wh{u^{[\ga{1};\dm]}}(\xi),\wh{u^{[\ga{2};\dm]}}(\xi),\dots,\wh{u^{[\ga{d_{\dm}};\dm]}}(\xi)\big].
	\ee
	
\end{itemize}

From \cite[Lemma 7]{dhacha}, it is not hard to deduce that
\be\label{DEF}
\FF_{r;\dm}(\xi)D_{u,\omega;\dm}(\xi)
\ol{\FF_{r;\dm}(\xi)}^{\tp}=d_{\dm} E_{u,\omega;\dm}(\dm^{\tp}\xi),\qquad \xi\in\dR, \omega\in \Omega_{\dm},
\ee	
where $\FF_{r;\dm}(\xi)$ is the following $(rd_{\dm})\times (rd_{\dm})$ matrix:
\be\label{Fourier}
\FF_{r;\dm}(\xi):=\left(e^{-i\ga{l}\cdot(\xi+2\pi\om{k})}I_r\right)_{1\leq l,k\leq d_\dm}.
\ee	
Thus we further deduce that
\be\label{DEF:2}P_{u;\dm}(\xi)=Q_{u;\dm}(\dm^\tp \xi) \FF_{r;\dm}(\xi),\ee
where $P_{u;\dm}(\xi):=\big[\wh{u}(\xi+2\pi\om{1}),
\wh{u}(\xi+2\pi\om{2}),\dots,\wh{u}(\xi+2\pi\om{d_{\dm}})\big]$ as
in \eqref{Pb}. 

Now let $\theta,\tilde{\theta},a,\tilde{a}\in\dlrs{0}{r}{r}$ and $b,\tilde{b}\in\dlrs{0}{s}{r}$ be finitely supported filters. Recall that $(\{a;b\},\{\tilde{a};\tilde{b}\})_{\Theta}$ (where $\Theta:=\theta^\star *\tilde{\theta}$) is a dual $\dm$-framelet filter bank if and only if \er{spectral} holds. Using \er{DEF:2} and $\FF_{r;\dm}\ol{\FF_{r;\dm}}^{\tp}=d_{\dm}I_{d_{\dm}r}$, it is straight forward to see that \er{spectral} is equivalent to
\be\label{dffb:coset:0}
\cN_{a,\tilde{a},\Theta}(\xi)=\ol{Q_{b;\dm}(\xi)}^\tp Q_{\tilde{b};\dm}(\xi),\ee
with
\be\label{dffb:coset}\cN_{a,\tilde{a},\Theta}(\xi):=
d_{\dm}^{-1}E_{\Theta,0;\dm}(\xi)-
\ol{Q_{a;\dm}(\xi)}^\tp \wh{\Theta}(\xi) Q_{\tilde{a};\dm}(\xi).
\ee
Therefore, constructing a dual framelet filter bank is equivalent to obtaining a matrix factorization as in \er{dffb:coset}. When the refinement masks $a$ and $\tilde{a}$ are given, all we have to do is to choose some suitable $\theta$ and $\tilde{\theta}$, and then factorize $\cN_{a,\tilde{a},\Theta}$ as in \er{dffb:coset}. Noting that the matrices $Q_{b;\dm}$ and $Q_{\tilde{b};\dm}$ give us all coset sequences of $b$ and $\tilde{b}$, we can finally reconstruct $b$ and $\tilde{b}$ via \er{f:coset}. It is worth mentioning that the approach of passing to coset sequences often appears in the literature of filter bank construction.\\

\bp[Proof of Theorem~\ref{thm:dfrt}] Since $\sr(a,\dm)=\tilde{m}$ and $\sr(\tilde{a},\dm)=m$ with matching filters $\vgu,\tvgu\in\dlrs{0}{1}{r}$ respectively, it is not hard to derive (see e.g. \cite[Proposition 3.2]{han03}) that 
$$\wh{\vgu}(\xi)\wh{\phi}(\xi)=\wh{\vgu}(0)\wh{\phi}(0)+\bo(|\xi|^{\tilde{m}}),\quad \wh{\tvgu}(\xi)\wh{\tphi}(\xi)=\wh{\tvgu}(0)\wh{\tphi}(0)+\bo(|\xi|^{m}),$$
as $\xi\to 0$. By proper scaling, WLOG we may assume that 
$$\wh{\vgu}(\xi)\wh{\phi}(\xi)=1+\bo(|\xi|^{\tilde{m}}),\quad \wh{\tvgu}(\xi)\wh{\tphi}(\xi)=1+\bo(|\xi|^{m}),$$
as $\xi\to 0$. Define $\wh{\Vgu_{\dn}}$ as in \er{vgu:special}. It follows from Theorem~\ref{thm:normalform:gen} that there exist strongly invertible filters $\theta,\tilde{\theta}\in\dlrs{0}{r}{r}$ such that 
$$\wh{\mrvgu}(\xi):=\wh{\vgu}(\xi)\wh{\theta}(\xi)^{-1}=r^{-1/2}\wh{\Vgu_{\dn}}(\xi)+\bo(\|\xi\|^{\tilde{m}}),$$ $$\wh{\mrphi}(\xi):=\wh{\theta}(\xi)\wh{\phi}(\xi)=r^{-1/2}\ol{\wh{\Vgu_{\dn}}(\xi)}^{\tp}+\bo(\|\xi\|^{m+\tilde{m}}),$$
$$\wh{\mrtvgu}(\xi):=\wh{\tvgu}(\xi)\wh{\tilde{\theta}}(\xi)^{-1}=r^{-1/2}\wh{\Vgu_{\dn}}(\xi)+\bo(\|\xi\|^{m}),$$ $$\wh{\mrtphi}(\xi):=\wh{\tilde{\theta}}(\xi)\wh{\tphi}(\xi)=r^{-1/2}\ol{\wh{\Vgu_{\dn}}(\xi)}^{\tp}+\bo(\|\xi\|^{m+\tilde{m}}),$$
as $\xi\to 0$. In particular, one can conclude that the following moment conditions hold as $\xi\to 0$:
\be\label{mr:dual}\wh{\mrvgu}(\xi)=C\ol{\wh{\mrtphi}(\xi)}^{\tp}+\bo(\|\xi\|^{\tilde{m}})=\wh{c}(\xi)\wh{\Vgu_{\dn}}(\xi)+\bo(\|\xi\|^{\tilde{m}}),\ee
\be\label{mrt:dual}\wh{\mrtvgu}(\xi)=\tilde{C}\ol{\wh{\mrphi}(\xi)}^{\tp}+\bo(\|\xi\|^m)=\wh{d}(\xi)\wh{\Vgu_{\dn}}(\xi)+\bo(\|\xi\|^m),\ee
\be\label{mrphi:mrtphi}\ol{\wh{\mrphi}(\xi)}^{\tp}\wh{\mrtphi}(\xi)=1+\bo(\|\xi\|^{m+\tilde{m}}),\ee	
for some $c,d\in\dlp{0}$ with $\wh{c}(0)\ne 0$ and $\wh{d}(0)\ne 0$, and some $C,\tilde{C}\in\C\setminus\{0\}$.

Define $\mra,\mrta\in\dlrs{0}{r}{r}$ as in \er{mr:filter} and \er{mrt:filter}. We have $\wh{\mrphi}(\dm^{\tp}\cdot)=\wh{\mra}\wh{\mrphi}$, and $\wh{\mrtphi}(\dm^{\tp}\cdot)=\wh{\mrta}\wh{\mrtphi}$. Furthermore, $\mra$ (resp. $\mrta$) has order $\tilde{m}$ (resp. $m$) sum rules with respect to $\dm$ with a matching filter $\mrvgu$ (resp. $\mrtvgu$).\\

Define $n:=\tilde{m}+m$. By Theorem~\ref{thm:normalform:gen}, there exists a strongly invertible $U\in\dlrs{0}{r}{r}$ such that 
$$\wh{\bpphi}(\xi):=\wh{U}(\xi)\wh{\mrphi}(\xi)=(1,0,\dots,0)^{\tp}+\bo(\|\xi\|^n),$$ $$\wh{\bpvgu}(\xi):=\wh{\mrvgu}(\xi)\wh{U}(\xi)^{-1}=(1,0,\dots,0)+\bo(\|\xi\|^{\tilde{m}}),$$
as $\xi\to 0$. Thus by letting $\wh{\bpa}:=\wh{U}(\dm^{\tp}\cdot)\wh{\mra}\wh{U}^{-1}$ , we see that $\bpa$ takes the ideal $(\tilde{m},n)$-normal form, that is, 
$$\wh{\bpa}(\xi)=\begin{bmatrix}\wh{\bpa_{1,1}}(\xi) &\wh{\bpa_{1,2}}(\xi)\\
\wh{\bpa_{2,1}}(\xi) &\wh{\bpa_{2,2}}(\xi)\end{bmatrix},$$
where $\wh{\bpa_{1,1}},\wh{\bpa_{1,2}},\wh{\bpa_{2,1}}$ and $\wh{\bpa_{2,2}}$ are $1\times 1, 1\times (r-1), (r-1)\times 1$ and $(r-1)\times (r-1)$ matrices of $2\pi\dZ$-periodic trigonometric polynomials such that
\begin{align*} &\wh{\bpa_{1,1}}(\xi)=1+\bo(\|\xi\|^n),\quad \wh{\bpa_{1,1}}(\xi+2\pi\omega)=\bo(\|\xi\|^{\tilde{m}}),\quad \xi\to 0,\quad\forall \omega\in\Omega_{\dm}\setminus\{0\},\\
&\wh{\bpa_{1,2}}(\xi+2\pi\omega)=\bo(\|\xi\|^{\tilde{m}}),\quad \xi\to 0,\quad\forall \omega\in\Omega_{\dm},\\ &\wh{\bpa_{2,1}}(\xi)=\bo(\|\xi\|^n),\quad \xi\to 0,
\end{align*}
as $\xi\to 0$, where $\Omega_{\dm}:=\{\omega_1,\dots,\om{d_{\dm}}\}$ is defined as \er{omega:dm}.\\

On the other hand,  we have
$$\wh{\bptvgu}(\xi):=\wh{\mrtvgu}(\xi)\ol{\wh{U}(\xi)}^{\tp}=\ol{\wh{\mrphi}(\xi)}^{\tp}\ol{\wh{U}(\xi)}^{\tp}+\bo(\|\xi\|^m)=(1,0,\dots,0)+\bo(\|\xi\|^m),$$
$$\wh{\bptphi}(\xi):=\ol{\wh{U}(\xi)}^{-\tp}\wh{\mrtphi}(\xi)=\ol{\wh{U}(\xi)}^{-\tp}\ol{\wh{\mrvgu}(\xi)}^{\tp}+\bo(\|\xi\|^{\tilde{m}})=(1,0,\dots,0)^{\tp}+\bo(\|\xi\|^{\tilde{m}}),$$
as $\xi\to 0$. Moreover, the condition \er{mrphi:mrtphi} implies that
$$\wh{\bptphi_1}(\xi)=1+\bo(\|\xi\|^n),\qquad \xi\to 0,$$
where $\bptphi_1$ is the first coordinate of $\bptphi$. Thus by letting $\wh{\bpta}:=\ol{\wh{U}(\dm^{\tp}\cdot)}^{-\tp}\wh{\mrta}\ol{\wh{U}}^{\tp}$, we see that $\wh{\bptphi}(\dm^{\tp}\cdot)=\wh{\bpta}\wh{\bptphi}$ and $\bpta$ has order $m$ sum rules with respect to $\dm$ with a matching filter $\bptvgu$.  Furthermore, we have
$$\wh{\bpta}(\xi)=\begin{bmatrix}\wh{\bpta_{1,1}}(\xi) &\wh{\bpta_{1,2}}(\xi)\\
\wh{\bpta_{2,1}}(\xi) &\wh{\bpta_{2,2}}(\xi)\end{bmatrix},$$ 
where $\wh{\bpta_{1,1}},\wh{\bpta_{1,2}},\wh{\bpta_{2,1}}$ and $\wh{\bpta_{2,2}}$ are $1\times 1, 1\times (r-1), (r-1)\times 1$ and $(r-1)\times (r-1)$ matrices of $2\pi\dZ$-periodic trigonometric polynomials such that
\begin{align*}	 &\wh{\bpta_{1,1}}(\xi)=1+\bo(\|\xi\|^n),\quad \wh{\bpta_{1,1}}(\xi+2\pi\omega)=\bo(\|\xi\|^m),\quad \xi\to 0,\quad\forall \omega\in\Omega_{\dm}\setminus\{0\},\\
&\wh{\bpta_{1,2}}(\xi+2\pi\omega)=\bo(\|\xi\|^m),\quad \xi\to 0,\quad\forall \omega\in\Omega_{\dm},\\ &\wh{\bpta_{2,1}}(\xi)=\bo(\|\xi\|^{\tilde{m}}),\quad \xi\to 0,
\end{align*}
as $\xi\to 0$.\\

For $j=1,\dots,d_{\dm}$, define
$$\wh{A_j}(\xi):=\td(\omega_j)I_r-\ol{\wh{\bpa}(\xi)}^{\tp}\wh{\bpta}(\xi+2\pi\omega_j),$$
where $\td$ is defined as \er{delta:seq}. We have
$$\wh{A_1}(\xi)=I_r-\ol{\wh{\bpa}(\xi)}^{\tp}\wh{\bpta}(\xi)=\begin{bmatrix}\wh{A_{1;1}}(\xi) & \wh{A_{1;2}}(\xi)\\
\wh{A_{1;3}}(\xi) &\wh{A_{1;4}}(\xi)\end{bmatrix},$$
where $\wh{A_{1;1}},\wh{A_{1;2}},\wh{A_{1;3}}$ and $\wh{A_{1;4}}$ are
$1\times 1, 1\times (r-1), (r-1)\times1$ and $(r-1)\times (r-1)$ matrices of $2\pi\dZ$-periodic trigonometric polynomials, satisfying the following moment conditions as $\xi\to 0$:
\begin{align*}
&\wh{A_{1;1}}(\xi)=1-\left(\ol{\wh{\bpa_{1,1}}(\xi)}\wh{\bpta_{1,1}}(\xi)+\ol{\wh{\bpa_{2,1}}(\xi)}^{\tp}\wh{\bpta_{2,1}}(\xi)\right)=\bo(\|\xi\|^n),\\
&\wh{A_{1;2}}(\xi)=-\ol{\wh{\bpa_{1,1}}(\xi)}\wh{\bpta_{1,2}}(\xi)-\ol{\wh{\bpa_{2,1}}(\xi)}^{\tp}\wh{\bpta_{2,2}}(\xi)=\bo(\|\xi\|^{m}),\\
&\wh{A_{1;3}}(\xi)=-\ol{\wh{\bpa_{1,2}}(\xi)}^{\tp}\wh{\bpta_{1,1}}(\xi)-\ol{\wh{\bpa_{2,2}}(\xi)}^{\tp}\wh{\bpta_{2,1}}(\xi)=\bo(\|\xi\|^{\tilde{m}}).\end{align*}
For $j=2,\dots,d_\dm$, we have
$$\wh{A_j}(\xi)=-\ol{\wh{\bpa}(\xi)}^{\tp}\wh{\bpta}(\xi+2\pi\om{j})=\begin{bmatrix}\wh{A_{j;1}}(\xi) & \wh{A_{j;2}}(\xi)\\
\wh{A_{j;3}}(\xi) &\wh{A_{j;4}}(\xi)\end{bmatrix},$$
where $\wh{A_{j;1}},\wh{A_{j;2}},\wh{A_{j;3}}$ and $\wh{A_{j;4}}$ are
$1\times 1, 1\times (r-1), (r-1)\times1$ and $(r-1)\times (r-1)$ matrices of $2\pi\dZ$-periodic trigonometric polynomials for each $j$, satisfying the following moment conditions as $\xi\to 0$:
\begin{align*}
&\wh{A_{j;1}}(\xi)=-\left(\ol{\wh{\bpa_{1,1}}(\xi)}\wh{\bpta_{1,1}}(\xi+2\pi \om{j})+\ol{\wh{\bpa_{2,1}}(\xi)}^{\tp}\wh{\bpta_{2,1}}(\xi+2\pi\om{j})\right)
=\bo(\|\xi\|^m),\\
&\wh{A_{j;1}}(\xi-2\pi \om{j})=-\left(\ol{\wh{\bpa_{1,1}}(\xi-2\pi \om{j})}\wh{\bpta_{1,1}}(\xi)+\ol{\wh{\bpa_{2,1}}(\xi-2\pi \om{j})}^{\tp}\wh{\bpta_{2,1}}(\xi)\right)=\bo(\|\xi\|^{\tilde{m}}),\\
&\wh{A_{j;2}}(\xi)=-\ol{\wh{\bpa_{1,1}}(\xi)}
\wh{\bpta_{1,2}}(\xi+2\pi\om{j})-\ol{\wh{\bpa_{2,1}}(\xi)}^{\tp}\wh{\bpta_{2,2}}(\xi+2\pi\om{j})=\bo(\|\xi\|^m),\\
&\wh{A_{j;3}}(\xi-2\pi\om{j})=-\ol{\wh{\bpa_{1,2}}(\xi-2\pi\om{j})}
\wh{\bpta_{1,1}}(\xi)-\ol{\wh{\bpa_{2,2}}(\xi-2\pi\om{j})}^{\tp}\wh{\bpta_{2,1}}(\xi)= \bo(\|\xi\|^{\tilde{m}}).
\end{align*}

For $\mu\in\dN_0$, define $\Delta_\mu\in \dlrs{0}{r}{r}$ via $\Delta_{\mu}:=\DG(\nabla^{\mu}\td, I_{r-1})$. From what we have done above, we conclude that
\be\label{bpaj:fac:dual}
\wh{A_j}(\xi)=\sum_{\alpha\in\dN_{0;m},\beta\in\dN_{0;{\tilde{m}}}}\ol{\wh{\Delta_\alpha}(\xi)}^{\tp}\wh{A_{j,\alpha,\beta}}(\xi)
\wh{\Delta_\beta}(\xi+2\pi\om{j}),\ee
for some $A_{j,\alpha,\beta}\in\dlrs{0}{r}{r}$ for all $\alpha\in\dN_{0;m},\beta\in\dN_{0;{\tilde{m}}}$ and all $j=1,\dots,d_{\dm}$.\\

Define $\cM_{\bpa,\bpta,I_r}$ as in \er{m:a:ta} with $a,\tilde{a}, \Theta$ being replaced by $\bpa,\bpta,I_r$ respectively, and recall that $D_{\mu,\omega;\dm}$ is defined as \er{Duni} for all $u\in\dlrs{0}{r}{r}$ and $\omega\in\Omega_{\dm}$. Note that $$\cM_{\bpa,\bpta,\tilde{U}}=\sum_{j=1}^{d_{\dm}}D_{A_j,\omega_j;\dm}=\sum_{\alpha\in\dN_{0;m},\beta\in\dN_{0;{\tilde{m}}}}\ol{D_{\Delta_\alpha,0;\dm}}^{\tp}D_{A_{j,\alpha,\beta},\omega_j;\dm}
D_{\Delta_\beta,0;\dm},$$
where the last identity follows from \er{bpaj:fac:dual}. \\

Define $\cN_{\bpa,\bpta,I_r}$ as in \er{dffb:coset} with $a,\tilde{a}$ and $\Theta$ being replaced by $\bpa,\bpta$ and $I_r$ respectively. Recall that $E_{\mu,\omega;\dm}$ is defined as \er{Euni} for all $u\in\dlrs{0}{r}{r}$ and $\omega\in\Omega_{\dm}$, and $\FF_{r;\dm}$ is defined as \er{Fourier}. It follows from \er{DEF:2} and $\FF_{r;\dm}\ol{\FF_{r;\dm}}^{\tp}=d_{\dm}I_{d_{\dm}r}$ that
\be\label{fac:cn}\begin{aligned}&\cN_{\bpa,\bpta,I_r}(\dm^{\tp}\xi)=d_{\dm}^{-2}\FF_{r;\dm}(\xi)\cM_{\bpa,\bpta,I_r}(\xi)\ol{\FF_{r;\dm}(\xi)}^{\tp}\\
	=&d_{\dm}^{-1}\sum_{j=1}^{d_{\dm}}
	\sum_{\alpha\in\dN_{0;m},\beta\in\dN_{0;\tilde{m}}}
	\ol{E_{\Delta_\alpha,0;\dm}(\dm^{\tp}\xi)}^{\tp}
	E_{A_{j,\alpha,\beta},\omega_j;\dm}(\dm^{\tp}\xi)
	E_{\Delta_\beta,0;\dm}(\dm^{\tp}\xi).
\end{aligned}\ee
By letting
$$E_{\alpha,\beta}(\xi):=d_{\dm}^{-1}\sum_{j=1}^{d_{\dm}}E_{A_{j,\alpha,\beta},\omega_j;\dm}(\xi),\qquad \xi\in\dR,\quad\alpha\in\dN_{0;m},\quad\beta\in\dN_{0;\tilde{m}}, $$
we have
\be\label{fac:cn:2}\cN_{\bpa,\bpta,I_r}(\xi)=
\sum_{\alpha\in\dN_{0;m},\beta\in\dN_{0;\tilde{m}}}
\ol{E_{\Delta_\alpha,0;\dm}(\xi)}^{\tp}
E_{\alpha,\beta}(\xi)
E_{\Delta_\beta,0;\dm}(\xi).\ee
For every $\alpha\in\dN_{0;m}$ and $\beta\in\dN_{0;\tilde{m}}$, choose $d_{\dm}r \times d_{\dm}r$ matrices $E_{\alpha,\beta,1}$ and $E_{\alpha,\beta,1}$ of $2\pi\dZ$-periodic trigonometric polynomials such that $E_{\alpha,\beta}=\ol{E_{\alpha,\beta,1}}^{\tp}E_{\alpha,\beta,2}$. Define $\bpb_{\alpha,\beta, k},\bptb_{\alpha,\beta,k}\in\dlrs{0}{1}{r}$ for $k=1,\dots, d_{\dm}r$ and all $\alpha\in\dN_{0;m},\beta\in\dN_{0;\tilde{m}}$ via
\begin{align}&\wh{\bpb_{\alpha,\beta}}(\xi):=\begin{bmatrix}\wh{\bpb_{\alpha,\beta,1}}(\xi)\\
\vdots\\
\wh{\bpb_{\alpha,\beta,d_{\dm}r}}(\xi)\end{bmatrix}:=E_{\alpha,\beta,1}(\dm^{\tp}\xi)\FF_{r;\dm}(\xi)\begin{bmatrix}\wh{\Delta_{\alpha}}(\xi)\\
\pmb{0}_{d_{\dm}(r-1)\times r}\end{bmatrix},\label{bab}\\
&\wh{\bptb_{\alpha,\beta}}(\xi):=\begin{bmatrix}\wh{\bptb_{\alpha,\beta,1}}(\xi)\\
\vdots\\
\wh{\bptb_{\alpha,\beta,d_{\dm}r}}(\xi)\end{bmatrix}:=E_{\alpha,\beta,2}(\dm^{\tp}\xi)\FF_{r;\dm}(\xi)\begin{bmatrix}\wh{\Delta_{\beta}}(\xi)\\
\pmb{0}_{d_{\dm}(r-1)\times r}\end{bmatrix},\label{btab}
\end{align}
where $\pmb{0}_{t\times q}$ denotes the $t\times q$ zero matrix. Recall that $P_{u;\dm}(\xi)=[\wh{u}(\xi+2\pi\omega_1),\dots,\wh{u}(\xi+2\pi\omega_{d_{\dm}})]$ as in \er{Pb} for all matrix-valued filter $u$. It is not hard to observe that
\be\label{Pbab}\begin{aligned}P_{\bpb_{\alpha,\beta};\dm}(\xi)&=E_{\alpha,\beta,1}(\dm^{\tp}\xi)\FF_{r;\dm}(\xi)D_{\Delta_{\alpha},0;\dm}(\xi)\\
	&=E_{\alpha,\beta,1}(\dm^{\tp}\xi)E_{\Delta_{\alpha},0;\dm}(\dm^{\tp}\xi)\ol{\FF_{r;\dm}(\xi)}^{\tp},\end{aligned}\ee 
where the last identity follows from \er{DEF} and $\FF_{r;\dm}\ol{\FF_{r;\dm}}^{\tp}=d_{\dm}I_{d_{\dm}r}$. Similarly,
\be\label{Pbtab}P_{\bptb_{\alpha,\beta};\dm}(\xi)=E_{\alpha,\beta,2}(\dm^{\tp}\xi)E_{\Delta_{\beta},0;\dm}(\dm^{\tp}\xi)\ol{\FF_{r;\dm}(\xi)}^{\tp}.\ee
It follows from \er{fac:cn}, \er{fac:cn:2}, \er{Pbab} and \er{Pbtab} that
\be\label{fac:cm}\begin{aligned}\cM_{\bpa,\bpta,I_r}(\xi)&=\ol{\FF_{r;\dm}(\xi)}^{\tp}\cN_{\bpa,\bpta,I_r}(\dm^{\tp}\xi)\FF_{r;\dm}(\xi)\\
	&=\sum_{\alpha\in\dN_{0;m},\beta\in\dN_{0;\tilde{m}}}\ol{P_{\bpb_{\alpha,\beta};\dm}(\xi)}^{\tp}P_{\bptb_{\alpha,\beta};\dm}(\xi).
\end{aligned}\ee
Define
\begin{align*}&\{\bpb_{\ell}:\ell=1,\dots,s\}:=\{\bpb_{\alpha,\beta}:\alpha\in\dN_{0;m},\quad\beta\in\dN_{0;\tilde{m}}\},\\
&\{\bptb_{\ell}:\ell=1,\dots,s\}:=\{\bptb_{\alpha,\beta}:\alpha\in\dN_{0;m},\quad\beta\in\dN_{0;\tilde{m}}\},\end{align*}
and let $\bpb:=[\bpb_1^{\tp},\dots,\bpb_s^{\tp}]^{\tp}, \bptb:=[\bptb_1^{\tp},\dots,\bptb_s^{\tp}]^{\tp}.$ We see that \er{fac:cm} becomes
$$\cM_{\bpa,\bpta,I_r}(\xi)=\ol{P_{\bpb;\dm}(\xi)}^{\tp}P_{\bptb;\dm}(\xi),$$
which is equivalent to say that $(\{\bpa;\bpb\},\{\bpta;\bptb\})_{I_r}$ is an OEP-based dual $\dm$-framelet filter bank satisfying
\be\label{oep:bp}\ol{\wh{\bpa}(\xi)}^{\tp}\wh{\bpta}(\xi+2\pi\omega)+\ol{\wh{\bpb}(\xi)}^{\tp}\wh{\bptb}(\xi+2\pi\omega)=\td(\omega)I_r,\qquad \xi\in\dR,\omega\in\Omega_{\dm}.\ee
Now define $\mrb,\mrtb,b,\tilde{b}\in\dlrs{0}{s}{r}$ via
$$\wh{\mrb}:=\wh{\bpb}\wh{U}^{-1},\quad \wh{\mrtb}:=\wh{\bptb}\ol{\wh{U}}^{\tp},\quad \wh{b}:=\wh{\mrb}\wh{\theta}^{-1},\quad \wh{\tilde{b}}:=\wh{\mrtb}\wh{\tilde{\theta}}^{-1}.$$
It follows from \er{oep:bp} that $(\{\mra;\mrb\},\{\mrta;\mrtb\})_{I_r}$ is an OEP-based dual $\dm$-framelet filter bank satisfying
$$\ol{\wh{\mra}(\xi)}^{\tp}\wh{\mrta}(\xi+2\pi\omega)+\ol{\wh{\mrb}(\xi)}^{\tp}\wh{\mrtb}(\xi+2\pi\omega)=\td(\omega)I_r,\qquad \xi\in\dR,\omega\in\Omega_{\dm},$$
and $(\{a;b\},\{\tilde{a};\tilde{b}\})_{\Theta}$ (where $\Theta:=\theta^{\star}*\tilde{\theta}$) is an OEP-based dual $\dm$-framelet filter bank satisfying \er{dffb}. By \er{mrt:dual} and \er{bab}, we have
\be\label{bpo:b}\begin{aligned}&\wh{\Vgu_{\dn}}(\xi)\ol{\wh{\mrb}(\xi)}^{\tp}=\ol{\wh{d}(\xi)^{-1}\wh{\mrphi}(\xi)}^{\tp}\ol{\wh{\mrb}(\xi)}^{\tp}+\bo(\|\xi\|^m)=\ol{\wh{d}(\xi)}^{-1}\ol{\wh{\bpphi}(\xi)}^{\tp}\ol{\wh{\bpb}(\xi)}^{\tp}+\bo(\|\xi\|^m)\\
	=&\ol{\wh{d}(\xi)}^{-1}(1,0,\dots,0)\ol{\wh{\bpb}(\xi)}^{\tp}+\bo(\|\xi\|^m)=\bo(\|\xi\|^m),\quad \xi\to 0,\end{aligned}\ee
where $d\in\dlp{0}$ with $\wh{d}(0)\ne 0$ is the same as in \er{mrt:dual}. 
Similarly, we deduce from \er{mr:dual} and \er{btab} that
\be\label{bpo:b:3}\wh{\Vgu_{\dn}}(\xi)\ol{\wh{\mrtb}(\xi)}^{\tp}=\bo(\|\xi\|^{m}),\qquad\xi\to 0.\ee
On the other hand, it follows immediately from \er{mrt:dual} and the refinement relation $\wh{\mrphi}(\dm^{\tp}\cdot)=\wh{\mra}\wh{\mrphi}$ that
$$\frac{\wh{d}(\xi)}{\wh{d}(\dm^{\tp}\xi)}\wh{\Vgu_{\dn}}(\xi)\ol{\wh{\mra}(\xi)}^{\tp}=\wh{\Vgu_{\dn}}(\dm^{\tp}\xi)+\bo(\|\xi\|^{\tilde{m}}),\qquad \xi\to 0.$$
Hence by Theorem~\ref{thm:bp}, we have $\bpo(\{\mra;\mrb\},\dm,\dn)=m=\sr(\mrta;\dm)$. This proves item (2).\\

Now define vector functions $\psi$ and $\psi$ as in \er{oep:mr} and \er{oep:mrt}. It follows from \er{mr:dual}, \er{mrt:dual}, \er{bpo:b} and \er{bpo:b:3} that $\vmo(\psi)=m$ and $\vmo(\tpsi)=\tilde{m}$. Further note that
$$\ol{\wh{\phi}(0)}^{\tp}\wh{\Theta}(0)\wh{\tphi}(0)=\ol{\wh{\mrphi}(0)}^{\tp}\wh{\mrtphi}(0)=1.$$
It follows from Theorem~\ref{thm:df} that $(\{\mrphi;\psi\},\{\mrtphi;\tpsi\})$ is a dual $\dm$-framelet in $\dLp{2}$. This proves item (3), and the proof is now complete.\ep

For the case $r=1$, it is too much to expect items (1) and (2) to hold, this is because a filter $\theta\in\dlp{0}$ is strongly invertible if and only if $\theta=c\td(\cdot-k)$ for some $c\in\C$ and $k\in\dZ$, and using such filters loses the advantage of OEP for increasing vanishing moments on framelet generators. Nevertheless, the matrix decomposition technique in the proof of Theorem~\ref{thm:dfrt} can be applied to deduce the following result for the case $r=1$.

\begin{cor}\label{cor:dft:r:1} Let $\dm$ be a $d\times d$ dilation matrix and let $\phi,\tphi\in \dLp{2}$ be compactly supported refinable functions satisfying $\wh{\phi}(\dm^{\tp} \xi)=\wh{a}(\xi)\wh{\phi}(\xi)$ and $\wh{\tphi}(\dm^{\tp} \xi)=\wh{\tilde{a}}(\xi)\wh{\tphi}(\xi)$, where $a,\tilde{a}\in\dlp{0}$ have order $\tilde{m}$ and $m$ sum rules with respect to $\dm$ with matching filters $\vgu,\tvgu\in \dlp{0}$, respectively. Suppose that $\wh{\vgu}(0)\wh{\phi}(0)=\wh{\tvgu}(0)\wh{\tphi}(0)=1$. Then there exist $b,\tilde{b}\in \dlrs{0}{s}{1}$ and $\theta,\tilde{\theta}\in \dlp{0}$ such that
	
	\begin{enumerate}
		\item  $(\{a;b\},\{\tilde{a};\tilde{b}\})_{\theta^\star*\tilde{\theta}}$ forms an OEP-based dual $\dm$-framelet filter bank.
		
		\item $(\{\mrphi;\psi\},\{\mrtphi;\tpsi\})$ is a compactly supported dual $\dm$-framelet in $\dLp{2}$, where $\mrphi,\psi,\mrtphi$ and $\tpsi$ are defined as in \er{oep:mr} and \er{oep:mrt}. Moreover, $\vmo(\psi)=m$ and $\vmo(\tpsi)=\tilde{m}$.
		
	\end{enumerate}
\end{cor}

\subsection{A Structural Characterization of OEP-based Balanced Dual Multiframelets}

The most important step to obtain an OEP-based dual multiframelet with high balancing orders is to determine a pair of balanced moment correction filters $\theta$ and $\tilde{\theta}$. In this subsection, we provide a complete criterion on balanced moment correction filters, and thus a comprehensive structural characterization of OEP-based dual multiframelets with high balancing orders is established.

The proof of Theorem~\ref{thm:dfrt} gives us some clues on the requirements of balanced moment correction filters. The following theorem states the sufficient conditions for obtaining an OEP-based dual framelet with all desired properties.

\begin{theorem}\label{thm:mcf}Let $\dm$ be a $d\times d$ dilation matrix and $r\ge 2$ be an integer. Let $\phi,\tilde{\phi}\in(\dLp{2})^r$ be compactly supported $\dm$ refinable vector functions associated with refinement masks $a,\tilde{a}\in\dlrs{0}{r}{r}$. Suppose that $\sr(a,\dm)=\tilde{m}$ and $\sr(\tilde{a},\dm)=m$ with matching filters $\vgu,\tvgu\in\dlrs{0}{1}{r}$ respectively such that $\wh{\vgu}(0)\wh{\phi}(0)\ne 0$ and $\wh{\tvgu}(0)\wh{\tphi}(0)\ne 0$. Let $\dn$ be a $d\times d$ integer matrix with $|\det(\dn)|=r$, and define $\wh{\Vgu_{\dn}}$ as in \er{vgu:special}.

	Let $\theta,\tilde{\theta}\in\dlrs{0}{r}{r}$ be strongly invertible finitely supported filters. Then $\theta$ and $\tilde{\theta}$ is a pair of $E_{\dn}$-balanced moment correction filters for the refinement masks $a$ and $\tilde{a}$ if 
	
	\begin{enumerate}
		\item[(i)] the moment conditions \er{mr:dual},\er{mrt:dual} and \er{mrphi:mrtphi} hold as $\xi\to 0$, for some $c,d\in\dlp{0}$ with $\wh{c}(0)\ne 0$ and $\wh{d}(0)\ne 0$, and some $C,\tilde{C}\in\C\setminus\{0\}$, where $\wh{\mrvgu}:=\wh{\vgu}\wh{\theta}^{-1},\wh{\mrphi}:=\wh{\theta}\wh{\phi},\wh{\mrtvgu}:=\wh{\tvgu}\wh{\tilde{\theta}}^{-1}$ and $\wh{\mrtphi}:=\wh{\tilde{\theta}}\wh{\tphi}$.
	\end{enumerate}
	
	Conversely, suppose $\theta$ and $\tilde{\theta}$ is a pair of $E_{\dn}$-balanced moment correction filters for the refinement masks $a$ and $\tilde{a}$, and suppose in addition the following statements hold:
	
	\begin{enumerate}
		\item[(ii)] $1$ is a simple eigenvalue of $\wh{a}(0)$ and $\wh{\mrta}(0)$. Moreover,  
		$$\lambda^{\alpha}I_r-\wh{a}(0),\quad I_r -\lambda^{\beta}\wh{a}(0),\quad I_r-\lambda^{\alpha}\wh{\tilde{a}}(0),\quad \lambda^{\beta}I_r -\wh{\tilde{a}}(0) $$ 
		are invertible matrices for all $\alpha, \beta\in\dN_0$ with $0<|\alpha|<\tilde{m}$ and $0<|\beta|<m$, where $\lambda:=(\lambda_1,\dots,\lambda_d)$ is the vector of the eigenvalues of $\dm$.

		\item[(iii)] $\wh{p}(\dm^{\tp}\xi)\wh{\Vgu_{\dn}}(\dm^{\tp}\xi)\wh{\mrta}(\xi)=\wh{p}(\xi)\wh{\Vgu_{\dn}}(\xi)+\bo(\|\xi\|^{m})$ as $\xi\to 0$ for some $p\in\dlp{0}$ with $\wh{p}(0)\neq 0$, where $\wh{\mrta}:=\wh{\tilde{\theta}}(\dm^{\tp}\cdot)\wh{\tilde{a}}\wh{\tilde{\theta}}^{-1}$ is defined as in \er{mrt:filter}.
		
		\item[(iv)] $\wh{q}(\xi)\wh{\mrta}(\xi)\ol{\wh{\Vgu_{\dn}}(\xi)}^{\tp}=\wh{q}(\dm^{\tp}\xi)\ol{\wh{\Vgu_{\dn}}(\dm^{\tp}\xi)}^{\tp}+\bo(\|\xi\|^{\tilde{m}})$ as $\xi\to 0$ for some $q\in\dlp{0}$ with $\wh{q}(0)\neq 0$.

	\end{enumerate}
	Then item (i) must hold.	
	
\end{theorem}

\bp  Suppose item (i) holds, then from the proof of Theorem~\ref{thm:dfrt}, one can immediately conclude that $\theta$ and $\tilde{\theta}$ is a pair of $E_{\dn}$-balanced moment correction filters for the refinement masks $a$ and $\tilde{a}$.

Conversely, suppose $\theta$ and $\tilde{\theta}$ is a pair of $E_{\dn}$-balanced moment correction filters for the refinement masks $a$ and $\tilde{a}$. Define $\mra,\mrta\in\dlrs{0}{r}{r}$ as in \er{mr:filter} and define
$\mrb,\mrtb\in\dlrs{0}{s}{r}$ as in \er{mrt:filter}. By item (2) of Theorem~\ref{thm:dfrt}, we have
\be\label{dffrt}\ol{\wh{\mra}(\xi)}^{\tp}\wh{\mrta}(\xi)+\ol{\wh{\mrb}(\xi)}^{\tp}\wh{\mrtb}(\xi)=I_r,\ee
and $\bpo(\{\mra;\mrb\},\dm,\dn)=m$. By Theorem~\ref{thm:bp}, we have
\be\label{bp:mr}\wh{\Vgu_{\dn}}(\xi)\ol{\wh{\mrb}(\xi)}^{\tp}=\bo(\|\xi\|^m),\quad \wh{\Vgu_{\dn}}(\xi)\ol{\wh{\mra}(\xi)}^{\tp}=\wh{\mrc}(\xi)\wh{\Vgu_{\dn}}(\dm^{\tp}\xi)+\bo(\|\xi\|^m),\quad \xi\to 0,\ee
for some $\mrc\in\dlp{0}$ with $\wh{\mrc}(0)\neq 0$.

Assume in addition that items (ii) - (iv) hold.

By left multiplying $\wh{\Vgu_{\dn}}$ on both sides of \er{dffrt} and using item (iv), we have
$$\wh{\Vgu_{\dn}}(\xi)=\wh{\mrc}(\xi)\wh{\Vgu_{\dn}}(\dm^{\tp}\xi)\wh{\mrta}(\xi)+\bo(\|\xi\|^m)=\wh{\mrc}(\xi)\frac{\wh{p}(\xi)}{\wh{p}(\dm^{\tp}\xi)}\wh{\Vgu_{\dn}}(\xi)+\bo(\|\xi\|^m),\quad\xi\to 0.$$
From the above relation we conclude that $\wh{\mrc}(0)=1$, and thus
\be\label{sr:mra}\wh{\mrd}(\dm^{\tp}\xi)\wh{\Vgu_{\dn}}(\dm^{\tp}\xi)\wh{\mrta}(\xi)=\wh{\mrd}(\xi)\wh{\Vgu_{\dn}}(\xi)+\bo(\|\xi\|^{m}),\quad \xi\to 0,\ee
where $\mrd\in\dlp{0}$ satisfies 
$$\wh{\mrd}(\xi)=\prod_{j=1}^\infty\wh{\mrc}((\dm^{\tp})^{-j}\xi)+\bo(\|\xi\|^m),\qquad\xi\to 0.$$
Moreover, it is easy to see from the second relation in \er{bp:mr} that 
\be\label{ref:mrta}\ol{\wh{\mrd}(\dm^{\tp}\xi)\wh{\Vgu_{\dn}}(\dm^{\tp}\xi)}^{\tp}=\wh{\mra}(\xi)\ol{\wh{\mrd}(\xi)\wh{\Vgu_{\dn}}(\xi)}^{\tp}+\bo(\|\xi\|^m),\quad\xi\to 0.\ee

We now apply the argument in the proof of \cite[Lemma 2.2]{han03} to prove that \er{mr:dual} and \er{mrt:dual} must hold.

Since $\mrta$ has $m$ sum rules with a matching filter $\mrtvgu$ with $\wh{\mrtvgu}:=\wh{\tvgu}\wh{\tilde{\theta}}^{-1}$, we have $\wh{\mrtvgu}(\dm^{\tp}\xi)\wh{\mrta}(\xi)=\wh{\mrtvgu}(\xi)+\bo(\|\xi\|^{m})$ as $\xi\to 0$. 
For any $n\in\N$, define $\otimes^nA:=A\otimes\dots\otimes A$ with $n$ copies of $A$. Recall that if $A,B,C$ and $E$ are matrices of sizes such that one can perform the matrix products $AC$ and $BE$, then we have $(A\otimes B)(C\otimes E)=(AC)\otimes(BE)$. Thus by induction, we have $(\otimes^n (AC))\otimes (BE)=[(\otimes^n A)\otimes B][(\otimes^n C)\otimes E]$.Furthermore, define the $1\times d$ vector of differential operators \be\label{diff:op}D:=(\partial_1,\dots,\partial_d),\text{ where }\partial_j:=\frac{\partial}{\partial \xi_j},\quad j=1,\dots, d.\ee 
It follows that
\be\label{tvgu:mom}\left[(\otimes^jD)\otimes\left(\wh{\mrtvgu}(\dm^{\tp}\cdot)\wh{\mrta}\right)\right](0)=[(\otimes^jD)\otimes\wh{\mrtvgu}](0),\qquad j=1,\dots,m-1.\ee
Rearranging \er{tvgu:mom} yields
\be\label{tvgu:mom:1}\left[(\otimes^jD)\otimes\left(\wh{\mrtvgu}(\dm^{\tp}\cdot)\wh{\mrta}(0)-\wh{\mrtvgu}\right)\right](0)=\left[(\otimes^jD)\otimes\left(\wh{\mrtvgu}(\dm^{\tp}\cdot)(\wh{\mrta}(0)-\wh{\mrta})\right)\right](0),\ee
for all $j=1,\dots,m-1.$ By the generalized product rule, we observe that the right hand of \er{tvgu:mom:1} only involves $\partial^\mu\wh{\mrtvgu}(0)$ with $|\mu|<j$. By calculation, we have
\be\label{tvgu:mom:2}\begin{aligned}&\left[(\otimes^jD)\otimes\left(\wh{\mrtvgu}(\dm^{\tp}\cdot)\wh{\mrta}(0)-\wh{\mrtvgu}\right)\right](0)	=\left([(\otimes^jD)\otimes \wh{\mrtvgu}](0)\right)[(\otimes^j\dm^{\tp})\otimes \wh{\mrta}(0)-I_{d^jr}],\end{aligned}\ee
for all $j\in\N$. Now by the condition in item (iii) on $\mrta$, the matrix $[(\otimes^j\dm^{\tp})\otimes \wh{\mrta}(0)-I_{d^jr}]$ is invertible for $j=1,\dots,m-1$. Moreover, it follows from \er{tvgu:mom:1} and \er{tvgu:mom:2} that up to a multiplicative constant, the values $\partial^\mu\wh{\mrtvgu}(0), |\mu|<m$ are uniquely determined via $\wh{\mrtvgu}(0)\wh{\mrta}(0)=\wh{\mrtvgu}(0)$ and
$$[(\otimes^jD)\otimes \wh{\mrtvgu}](0)=\left(\left[(\otimes^jD)\otimes\left(\wh{\mrtvgu}(\dm^{\tp}\cdot)(\wh{\mrta}(0)-\wh{\mrta})\right)\right](0)\right)[(\otimes^j\dm^{\tp})\otimes \wh{\mrta}(0)-I_{d^jr}]^{-1},$$
for all $j=1,\dots,m-1$.

Next, note that the refinement relation $\wh{\mrphi}(\dm^{\tp}\cdot)=\wh{\mra}\wh{\mrphi}$ (where $\wh{\mrphi}=\wh{\theta}\wh{\phi}$) holds. This implies that
\be\label{mrphi:mom}[(\otimes^j D)\otimes \wh{\mrphi}(\dm^{\tp}\cdot)](0)=[(\otimes^jD)\otimes(\wh{\mra}\wh{\mrphi})](0),\qquad j\in\N.\ee
Rearranging \er{mrphi:mom} yields
\be\label{mrphi:mom:1}\left[(\otimes^j D)\otimes \left(\wh{\mrphi}(\dm^{\tp}\cdot)^{\tp}-\wh{\mrphi}^{\tp}\wh{\mra}(0)^{\tp}\right)\right](0)=\left[(\otimes^jD)\otimes\left(\wh{\mrphi}^{\tp}(\wh{\mra}^{\tp}-\wh{\mra}(0)^{\tp})\right)\right](0),\ee
for all $ j\in\N$. Note that the right hand side of \er{mrphi:mom:1} only involves $\partial^\mu\wh{\mrphi}(0)$ with $|\mu|<j$. Furthermore, direct calculation yields
\be\label{mrphi:mom:2}\begin{aligned}&\left[(\otimes^j D)\otimes \left(\wh{\mrphi}(\dm^{\tp}\cdot)^{\tp}-\wh{\mrphi}^{\tp}\wh{\mra}(0)^{\tp}\right)\right](0)\\
	=&\left([(\otimes^jD)\otimes\wh{\mrphi}^{\tp}](0)\right)[(\otimes^j\dm^{\tp})\otimes I_r-(\otimes^jI_d)\otimes \wh{\mra}(0)],\end{aligned}\ee
for all $j\in\N$. Now by the condition in item (iii) on $\mra$, the matrix $[(\otimes^j\dm^{\tp})\otimes I_r-(\otimes^jI_d)\otimes \wh{\mra}(0)]$ is invertible for $j=1,\dots,m-1$. Moreover, it follows from \er{mrphi:mom:1} and \er{mrphi:mom:2} that up to a multiplicative constant, the values $\partial^\mu\wh{\mrphi}(0), |\mu|<m$ are uniquely determined via $\wh{\mrphi}(0)=\wh{\mra}(0)\wh{\mrphi}(0)$ and
$$\begin{aligned}&[(\otimes^jD)\otimes\wh{\mrphi}^{\tp}](0)\\
=&\left(\left[(\otimes^jD)\otimes\left(\wh{\mrphi}^{\tp}(\wh{\mra}^{\tp}-\wh{\mra}(0)^{\tp})\right)\right](0)\right)[(\otimes^j\dm^{\tp})\otimes I_r-(\otimes^jI_d)\otimes \wh{\mra}(0)]^{-1},\end{aligned}$$
for all $j=1,\dots,m-1.$

Consequently, by the above analysis and using \er{sr:mra} and \er{ref:mrta}, we conclude that \er{mrt:dual} holds for some $\tilde{C}\in\C\setminus\{0\}$, with $d\in\dlp{0}$ being a non-zero scalar multiple of $\mrd$.\\

On the other hand, the condition on $\mrta$ in item (ii) and item (iv) together yield
\be\label{mrtphi:q}\wh{\mrtphi}(\xi)=K\wh{q}(\xi)\ol{\wh{\Vgu_{\dn}(\xi)}}^{\tp}+\bo(\|\xi\|^{\tilde{m}}),\qquad \xi\to 0,\ee
for some non-zero constant $K$. As $\theta$ and $\tilde{\theta}$ is a pair of $E_{\dn}$-balanced moment correction filters for the refinement masks $a$ and $\tilde{a}$, then in particular item (3) of Theorem~\ref{thm:dfrt} holds. Then $\vmo(\tpsi)=\tilde{m}$ and \er{mrtphi:q} imply that $\wh{\mrtb}(\xi)\ol{\wh{\Vgu_{\dn}}(\xi)}^{\tp}=\bo(\|\xi\|^{\tilde{m}})$ as $\xi \to 0$. Now right multiplying $\wh{q}\ol{\wh{\Vgu_{\dn}}}^{\tp}$ to both sides of \er{dffrt} yields
\be\label{ref:mra}\ol{\wh{q}(\dm^{\tp}\xi)}\wh{\Vgu_{\dn}}(\dm^{\tp}\xi)\wh{\mra}(\xi)=\ol{\wh{q}(\xi)}\wh{\Vgu_{\dn}}(\xi)+\bo(\|\xi\|^{\tilde{m}}),\quad \xi\to 0.\ee
Since $\mra$ has $\tilde{m}$ sum rules with a matching filter $\wh{\mrvgu}:=\wh{\vgu}\wh{\theta}^{-1}$, we have $\wh{\mrvgu}(\dm^{\tp}\xi)\wh{\mra}(\xi)=\wh{\mrvgu}(\xi)+\bo(\|\xi\|^{\tilde{m}})$ as $\xi\to 0$. Moreover, $\mrta$ satisfies the refinement equation $\wh{\mrtphi}(\dm^{\tp}\cdot)=\wh{\mrta}\wh{\mrtphi}$. By the condition in item (iii) on $\mra$, we conclude from \er{mrtphi:q} and \er{ref:mra} that \er{mr:dual} must hold for some $C\in\C\setminus\{0\}$, with $c\in\dlp{0}$ being a non-zero scalar multiple of $q$.\\

Finally, by left multiplying $\ol{\wh{\mrphi}}^{\tp}$ and right multiplying $\wh{\mrtphi}$ (where $\wh{\mrtphi}:=\wh{\tilde{\theta}}\wh{\tphi}$) to \er{dffrt}, we have
$$\ol{\wh{\mrphi}(\dm^{\tp}\xi)}^{\tp}\wh{\mrtphi}(\dm^{\tp}\xi)=\ol{\wh{\mrphi}(\xi)}^{\tp}\wh{\mrtphi}(\xi)+\bo(\|\xi\|^{\tilde{m}+m}),\quad\xi\to 0,$$
and \er{mrphi:mrtphi} follows from the above identity. This proves item (i), and the proof is now complete.\ep

\section{An Example of Balanced Moment Correction Filters}\label{sec:exmp}

\begin{exmp}Consider the compactly supported $M_{\sqrt{2}}$-refinable vector function $\phi=(\phi_1,\phi_2)^{\tp}$ whose refinement mask $a\in (l_0(\Z^2))^{2\times 2}$ is given by 
	$$\wh{a}(\xi)=\dfrac{1}{4}\begin{bmatrix}2 & 1+e^{i\xi_1}+e^{i\xi_2}+e^{i(\xi_1+\xi_2)}\\
	2e^{-i\xi_1} &0\end{bmatrix}\quad\text{and}\quad M_{\sqrt{2}}=\begin{bmatrix}1 & 1\\
	1& -1\end{bmatrix}.$$
The filter $a$ has sum rule of order $2$ with respect to $M_{\sqrt{2}}$, with a matching filter $\vgu\in (l_0(\Z^2))^{1\times 2}$ satisfying
$$\wh{\vgu}(\xi)=\left(1,1+\frac{i}{2}(\xi_1+\xi_2)\right)+\bo(\|\xi\|^2),\quad \xi=(\xi_1,\xi_2)\to(0,0).$$
The refinable vector function $\phi$ satisfies the moment condition
$$\wh{\phi}(\xi)=\frac{1}{360}\begin{bmatrix}240-12(\xi_1^2+\xi_2^2)\\[0.3cm]
120-60i(\xi_1+\xi_2)-6(3\xi_1^2+5\xi_1\xi_2+3\xi_2^2+i(4\xi_1^3+9\xi_1^2\xi_2+9\xi_1\xi_2^2+4\xi_2^3)) \end{bmatrix}+\bo(\|\xi\|^4)$$
as $\xi=(\xi_1,\xi_2)\to(0,0).$\\

\vspace{0.2cm}

For the dual refinable vector function, we take $\tphi=(\tphi_1,\tphi_2)^\tp$ whose refinement mask $\tilde{a}\in(l_0(\Z^2))^{2\times 2}$ is given by 
$$\wh{\tilde{a}}(\xi)=\dfrac{1}{16}\begin{bmatrix}8 & 3+e^{-i\xi_1}+e^{-i\xi_2}+3e^{-i(\xi_1+\xi_2)}\\
8e^{i\xi_1} & 1+3e^{i\xi_1}+3e^{-i\xi_2}+e^{i(\xi_1-\xi_2)}\end{bmatrix}.$$
The filter $\tilde{a}$ has sum rule of order $2$ with respect to $M_{\sqrt{2}}$, with a matching filter $\tvgu\in (l_0(\Z^2))^{1\times 2}$ satisfying
$$\wh{\tvgu}(\xi)=\left(1,1-\frac{i}{2}(\xi_1+\xi_2)\right)+\bo(\|\xi\|^2),\quad \xi=(\xi_1,\xi_2)\to(0,0).$$
The refinable vector function $\tphi$ satisfies the moment condition
$$\wh{\tphi}(\xi)=\frac{1}{192}\begin{bmatrix}96-2(4\xi_1^2+\xi_1\xi_2+2\xi_2^2)\\[0.3cm]
96+48i(\xi_1+\xi_2)-2(10\xi_1^2+13\xi_1\xi_2+8\xi_2^2-i(6\xi_1^3+11\xi_1^2\xi_2+9\xi_1\xi_2^2+4\xi_2^3)) \end{bmatrix}+\bo(\|\xi\|^4)$$
as $\xi=(\xi_1,\xi_2)\to(0,0).$\\

\vspace{0.2cm}

Let $\dn=M_{\sqrt{2}}$, we have $\wh{\Vgu_\dn}(\xi)=(1,e^{i(\xi_1+\xi_2)/2})$. One can construct a pair of $E_{\dn}$-balanced moment correction filters $\theta,\tilde{\theta}\in(l_0(\Z^2))^{2\times 2}$ given by
$$\wh{\theta}(\xi)=\begin{bmatrix}p_1(\xi) &p_2(\xi)\\
p_3(\xi) &p_4(\xi)\end{bmatrix},\quad \wh{\tilde{\theta}}(\xi)=\begin{bmatrix}q_1(\xi) &q_2(\xi)\\
q_3(\xi) &q_4(\xi)\end{bmatrix}\begin{bmatrix}-d_3(\xi) &d_2(\xi)\\
d_2(\xi) &-d_1(\xi)\end{bmatrix}\begin{bmatrix}c_1(\xi) &c_2(\xi)\\
c_2(\xi) &c_3(\xi)\end{bmatrix}\begin{bmatrix}s_1(\xi) &s_2(\xi)\\
s_3(\xi) &s_4(\xi)\end{bmatrix},$$
where $p_j,q_j,s_j, j=1,2,3,4$ and $d_k,c_k, k=1,2,3$ are $2\pi\Z^2$-periodic trigonometric polynomials given by the following:
\begin{align*}
	&p_1(\xi):=\frac{\sqrt{2}}{272}
	\left(542001-3225e^{-2i\xi_1}-7740e^{-i(\xi_1+\xi_2)}
	-265735e^{-i\xi_1}+12267e^{_i\xi_2}
	-4522e^{i(\xi_1-\xi_2)}-273258e^{i\xi_1}\right),\\
	&p_2(\xi):=-\frac{\sqrt{2}}{17}e^{-i(\xi_1+\xi_2)}\left(645e^{-i\xi_1}-646\right),\\
	&p_3(\xi):=\frac{\sqrt{2}}{3264}e^{i(\xi_1+\xi_2)}\left[(7740e^{-2i\xi_1}-12267e^{-i\xi_1}+4522)e^{-2i\xi_2}
	-1075e^{-2i\xi_1}-89655e^{-i\xi_1}+90873\right.\\
	&\qquad\qquad\qquad \qquad\qquad \left.+(3225e^{-3i\xi_1}+265735e^{-2i\xi_1}-544581e^{-i\xi_1}+274763)e^{-i\xi_2}
	\right],\\
	&p_4(\xi):=\frac{\sqrt{2}}{204}e^{-i\xi_1}(645e^{-i(\xi_1+\xi_2)}-646e^{-i\xi_2}-215).
	\end{align*}

\begin{align*}
q_1(\xi):=&\frac{\sqrt{2}}{2350080}
\left(152782870420
-171803098950e^{-i\xi_1}-24271223880e^{-i\xi_2}-473689039775e^{i\xi_1}\right.\\
&\left.+75422352630e^{i\xi_2}+2461284145e^{-i(\xi_1+\xi_2)}-3502342314e^{i(\xi_1+\xi_2)}+7992104616e^{i(\xi_1-\xi_2)}\right.\\
&\left.-4616116784e^{i(\xi_2-\xi_1)}+7473040038e^{-2i\xi_1}-9018774e^{-2i\xi_2}-8342459019e^{-3i\xi_1}+1511e^{-3i\xi_2}\right.\\
&\left.+8967362e^{-i(\xi_1+2\xi_2)}-8333493168e^{-i(2\xi_1+\xi_2)}+576060810e^{-i(2\xi_1-\xi_2)}\right),\\
q_2(\xi):=&\frac{\sqrt{2}}{8}\left(-6+e^{-i\xi_1}+e^{-i\xi_2}\right),\\
q_3(\xi):=&\frac{\sqrt{2}}{2350080}
\left(-41819386604
+25602071226e^{-i\xi_1}-24343301544e^{-i\xi_2}+16567932951e^{i\xi_1}\right.\\
&\left.-1735780129e^{i\xi_2}+24684568266e^{-i(\xi_1+\xi_2)}+1167447438e^{i(\xi_1+\xi_2)}+7992104616e^{i(\xi_1-\xi_2)}\right.\\
&\left.-7630304e^{i(\xi_2-\xi_1)}+7990728234e^{-2i\xi_1}-9006686e^{-2i\xi_2}-8342459019e^{-3i\xi_1}+1511e^{-3i\xi_2}\right.\\
&\left.+8967362e^{-i(\xi_1+2\xi_2)}-8333493168e^{-i(2\xi_1+\xi_2)}+576060810e^{-i(2\xi_1-\xi_2)}\right),\\
q_4(\xi):=&\frac{\sqrt{2}}{8}\left(2+e^{-i\xi_1}+e^{-i\xi_2}\right).
\end{align*}

\begin{align*}
d_1(\xi):=&\frac{1}{1152}\left(846168+964186e^{i(\xi_1+\xi_2)}-581936e^{i\xi_1}+648973e^{i(\xi_1-\xi_2)}-65e^{2i\xi_2}-1876500e^{i\xi_2}\right.\\
&\left.-1253380e^{-i\xi_2}+912657e^{i(\xi_2-\xi_1)}+576e^{-i\xi_1}+604428e^{-i(\xi_1+\xi_2)}+263955e^{-2i\xi_1}\right),\\
d_2(\xi):=&(1-d_1(\xi)/d_1(0))^2,\quad d_3(\xi):=\frac{d_2(\xi)^2-1}{d_1(\xi)}.
\end{align*}

\begin{align*}
c_1(\xi):=&\frac{1}{96}\left(64e^{i(\xi_1+2\xi_2)}-208e^{i(\xi_1+\xi_2)}+203e^{i\xi_1}-83e^{i(\xi_1-\xi_2)}-62e^{2i\xi_2}+199e^{i\xi_2}-58+109e^{-2i\xi_2}\right.\\
&\left.-17e^{i(\xi_2-\xi_1)}+49e^{-i\xi_1}-4e^{_2i\xi_1}\right),\\
c_2(\xi):=&(1-c_1(\xi)/c_1(0))^2,\quad c_3(\xi):=\frac{c_2(\xi)^2-1}{c_1(\xi)}.
\end{align*}

\begin{align*}
s_1(\xi):=&\frac{1}{8}(6-e^{-i\xi_1}-e^{-i\xi_2}),\quad s_2(\xi):=\frac{1}{8}(2+e^{-i\xi_1}+e^{-i\xi_2}),\\
s_3(\xi):=&\frac{1}{256}\left(-42e^{-i\xi_1}-42e^{-i\xi_2}-2e^{i(\xi_1-\xi_2)}+6e^{i(\xi_1+\xi_2)}+160+11e^{i\xi_2}+6e^{-i(\xi_1+\xi_2)}-2e^{i(\xi_2-\xi_1)}\right.\\
&\left.+11e^{i\xi_1}+12e^{-2i\xi_1}+12e^{-2i\xi_2}-e^{-3i\xi_1}-e^{-3i\xi_2}\right),\\
s_4(\xi):=&\frac{1}{256}\left(-6e^{-i\xi_1}-6e^{-i\xi_2}+2e^{i(\xi_1-\xi_2)}+2e^{i(\xi_1+\xi_2)}-112+5e^{i\xi_2}-14e^{-i(\xi_1+\xi_2)}+2e^{i(\xi_2-\xi_1)}\right.\\
&\left.+5e^{i\xi_1}-4e^{-2i\xi_1}-4e^{-2i\xi_2}+e^{-3i\xi_1}+e^{-3i\xi_2}\right).
\end{align*}

Both $\theta$ and $\tilde{\theta}$ are strongly invertible. Define $\wh{\mrvgu}:=\wh{\vgu}\wh{\theta}^{-1}$, $\wh{\mrphi}:=\wh{\theta}\wh{\phi}$, $\wh{\mrtvgu}:=\wh{\tvgu}\wh{\tilde{\theta}}^{-1}$, $\wh{\mrtphi}:=\wh{\tilde{\theta}}\wh{\tphi}$, we have
\begin{align*}\wh{\mrvgu}(\xi)&=\ol{\wh{\mrphi}(\xi)}^{\tp}+\bo(\|\xi\|^2)=\ol{\wh{\mrtphi}(\xi)}^{\tp}+\bo(\|\xi\|^2)=\wh{\mrtvgu}(\xi)+\bo(\|\xi\|^2)\\
&=\wh{c}(\xi)\wh{\Vgu_{\dn}}+\bo(\|\xi\|^2),\quad \xi\to(0,0),\end{align*}
where $\wh{c}$ is some $2\pi\Z^2$-periodic trigonometric polynomial satisfying
$$\wh{c}(\xi)=-\frac{\sqrt{2}}{2}+\frac{143\sqrt{2}}{8}i\xi_1-\frac{\sqrt{2}}{24}i\xi_2+\bo(\|\xi\|^2),\quad \xi\to (0,0),$$
and we have $\ol{\wh{\mrtphi}(\xi)}^\tp\wh{\mrphi}(\xi)=1+\bo(\|\xi\|^4)$ as $\xi\to (0,0)$.
\end{exmp}

\section*{Acknowledgements}
Part of the results of the paper was established under the guidance of my Ph.D. supervisor Dr. Bin Han. I wish to express my sincerest gratitude to him. Also, I would like to thank the editors of the journal for handling the manuscript, as well as the anonymous reviewers for carefully reviewing the contents.

\end{document}